\documentclass[11pt]{amsart}
\usepackage{ConvergenceEuler}

\author{Thomas O. Gallou\"et}
\address{CMLS, \'Ecole polytechnique, CNRS, Universit\'e Paris-Saclay, 91128 Palaiseau Cedex, France. E-mail: {\tt
 thomas.gallouet@polytechnique.edu}}
\thanks{The first author is supported by the ANR grant ISOTACE}

\author{Quentin M\'erigot}
\address{CNRS / Universit\'e Paris-Dauphine, Paris, France. E-mail: {\tt
 quentin.merigot@univ-dauphine.fr}}

\title[A Lagrangian scheme for the incompressible Euler equation]{A
  Lagrangian scheme for the incompressible Euler equation using
  optimal transport}


\begin{document}

\begin{abstract}
We approximate the regular solutions of the incompressible Euler
equation by the solution of ODEs on finite-dimensional spaces. Our
approach combines Arnold's interpretation of the solution of Euler's
equation for incompressible and inviscid fluids as geodesics in the
space of measure-preserving diffeomorphisms, and an extrinsic
approximation of the equations of geodesics due to Brenier.  Using
recently developed semi-discrete optimal transport solvers, this
approach yields numerical scheme able to handle problems of realistic
size in 2D. Our purpose in this article is to establish the convergence
of these scheme towards regular solutions of the incompressible Euler
equation, and to provide numerical experiments on a few simple
testcases in 2D.
\end{abstract}

\maketitle

\tableofcontents

\section{Introduction}

In this paper we investigate a discretization of Euler's equation for
incompressible and inviscid fluids in a domain $\Omega\subseteq \R^d$
with Neumann boundary conditions:
\begin{equation}
\begin{cases}
\partial_t v(t,x)  + \left( v(t,x) \cdot \nabla \right)v(t,x)=-\nabla p(t,x), &\hbox{ for } t \in[0,T],\;x\in  \Omega  \, ,  \medskip\\
\div (v(t,x)) = 0  &\hbox{ for } t \in[0,T],\;x\in  \Omega  \, ,  \medskip\\
v(t,x) \cdot  n = 0  &\hbox{ for } t \in[0,T],\;x\in  \partial \Omega  \, ,  \medskip\\
v(0,x) = v_0.
\end{cases} \label{eq:EI}
\end{equation}
As noticed by Arnold \cite{arnold1966geometrie}, in Lagrangian
coordinates, Euler's equation can be interpreted as the equation of
geodesics in the infinite-dimensional group of measure-preserving
diffeomorphisms of $\Omega$. To see this, we consider the flow map
$\phi: [0,T]\times \Omega \to \Omega$ induced by the vector field $v$,
that is:
\begin{equation}
\begin{cases}
\frac {d}{dt}  \phi(t,x)= v\left(t,\phi(t,x)\right) &\hbox{ for } t \in[0,T],\;x\in  \Omega  \, ,  \medskip\\
\phi(0,\cdot)=\id, \\
\partial_t \phi(0,\cdot)=v_0.
\end{cases} \label{eq:EIEtoL}
\end{equation}
Using the incompressibility constraint $\div(v(t,x))= 0$ and the
initial condition $\phi(0) = \id$, one can check that $\phi(t, \cdot)$
belongs to the set of volume preserving maps $\mS$, defined by
$$\mS =\left \{ s \in \LL^2(\Omega,\Rsp^d) \mid s_\# \Leb = \Leb
\right \},$$ where $\Leb$ is the restriction of the Lebesgue measure
to the domain $\Omega$ and where the pushforward measure $s_\#\Leb$ is
defined by the formula $s_\#\Leb(A) = \Leb(s^{-1}(A))$ for every
measurable subset $A$ of $\Omega$.  Euler's equation \eqref{eq:EI} can
therefore be reformulated as
\begin{equation}
\begin{cases}
\frac {d^2}{dt^2} \phi(t)= -\nabla  p(t,\phi(t,x)) &\hbox{ for } t \in[0,T],\;x\in  \Omega  \, ,  \medskip\\
\phi(t,\cdot) \in \mS &\hbox{ for } t \in[0,T], \medskip\\
\phi(0,\cdot)= \id , \medskip\\
\partial_t \phi(0,\cdot)=v_0.
\end{cases} \label{eq:EIL}
\end{equation}
This equation can be formally interpreted as the equation of geodesics
in $\mS$ as follows. First, note that the condition $\phi(t,\cdot) \in
\mS$ in \eqref{eq:EI} encodes the infinitesimal conditions $\div v(t,
\cdot) = 0$ and $v(t,x)\cdot n(x) = 0$ in \eqref{eq:EIL}. This
suggests that the tangent plane to $\mS$ at a point $\phi \in \mS$
should be the set $\{ v \circ \phi \mid v \in \mdiv\}$, where $\mdiv$
denotes the set of divergence-free vector fields
$$\mdiv =\left \{ v \in \LL^2(\Omega, \R^d) \mid \div(v)=0, \, v\cdot
n =0 \right \}.$$ In addition, by the Helmoltz-Hodge decomposition,
the orthogonal to $\mdiv$ in $\LL^2(\Omega,\Rsp^d)$ is the space of
gradients of functions in $\H^1_0(\Omega)$. Therefore the evolution
equation in \eqref{eq:EIL} expresses that the acceleration of $\phi$
should be orthogonal to the tangent plane to $\mS$ at $\phi$, or in
other words that $t \mapsto \phi(t,\cdot)$ should be a geodesic of
$\mS$.  Note however that a solution to \eqref{eq:EIL} does not need
to be a \emph{minimizing} geodesic between $\phi(0,\cdot)$ and
$\phi(T,\cdot)$.  The problem of finding minimizing geodesic on $\mS$
between two measure preserving maps, amounts to solving equations \eqref{eq:EIL}, where the initial condition $\partial_t \phi(0,\cdot)=v_0$ is replaced by a prescribed coupling between the position of particles at initial and final times.
It leads to generalized and
non-deterministic solutions introduced Brenier \cite{Brenier:1985bz},
where particles are allowed to split and cross. Shnirel'man showed
that this phenomena can happen even when the measure-preserving maps
$\phi(0,\cdot)$ and $\phi(T,\cdot)$ are diffeomorphisms of $\Omega$
\cite{Shnirelcprimeman:1994ef}.

Our discretization of Euler's equations \eqref{eq:EI} relies on
Arnold's interpretation as the equation of geodesics and exploit the
extrinsic view given by the embedding of the set of measure preserving
maps $\mS$ in the Hilbert space $\mM = \LL^2(\Omega,\Rsp^d)$. In our
discretization the measure-preserving property is enforced through a
penalization term involving the squared distance to the set of
measure-preserving maps $\mS$, as in \cite{brenier2000derivation}. The numerical implementation of this
idea relies on Brenier's polar factorization theorem to compute the
squared distance to $\mS$ and on recently developed numerical solvers
for optimal transport problems invoving a probability measure with
density and a finitely-supported probability measure
\cite{aurenhammer1998minkowski,merigot2011multiscale,de2012blue,levy2014numerical}.
This combination of ideas presented above has already been used to
compute numerically minimizing geodesics between measure-preserving
maps in \cite{mermir2015minimal}, allowing the recovery of
non-deterministic solutions predicted by Schnirel'man and Brenier. The
object of this article is to determine whether this strategy can be
used to construct a Lagrangian discretization for the more classical
Cauchy problem for the Euler's equation \eqref{eq:EI}, which is able
to recover regular solutions to Euler's equation, both theoretically and experimentally. 

\subsection*{Discretization in space: approximate geodesics}
The construction of approximate geodesics presented here is strongly
inspired by a particle scheme introduced by
Brenier~\cite{brenier2000derivation}, in which the space of
measure-preserving maps $\mS$ was approximated by the space of
permutations of a fixed tessellation of the domain $\Omega$. To
construct our numerical approximation we first approach the Hilbert
space $\mM = \LL^2(\Omega,\Rsp^d)$ with finite dimensional
subspaces. Let $N$ be an integer and let $P_N$ be a tessellation partition up to negligible set of
$\Omega$ into $N$ subsets $(\omega_i)_{1\leq i\leq N}$ satisfying
$$
\left\{
\begin{aligned}
  &\forall i \in\{1,\hdots, N\},~~ \Leb(\omega_i) = \frac{1}{N} \Leb(\Omega)\\
  & h_N := \max_{1 \leq i\leq N} \diam(\omega_i) \leq \frac{C}{N^{1/d}}
\end{aligned}
\right.$$ where $C$ is independent of
$N$. We consider $\mM_N$ the space of functions from $\Omega$ to
$\Rsp^d$ which are constant on each of the subdomains $(\omega_i)$. To
construct our approximate geodesics, we consider the squared distance
to the set $\mS\subseteq \mM$ of measure-preserving maps:
$$ d^2_{\mS}: m \in \mM_n \mapsto \min_{s\in \mS} \|{m -
  s}\|^2_{\mM}.$$ The approximate geodesic model is described by the
equations
\begin{equation}
\begin{cases}
\ddot m(t) + \frac{\nabla \dd^2_{\mS}(m(t))}{2\epsilon^2} =0, &\hbox{ for }t \in [0,T]  \, ,  \medskip\\
(m(0),\dot{m}(0)) \in  \mM_N^2
\end{cases} \label{eq:appgeogintro}
\end{equation}
which is the system associated to the Hamiltonian
\begin{equation} \label{Hamilto}
H(m,v)=\frac{1}{2}\nr{v}_\mM^2+\frac{d^2_{\mS}(m)}{2\epsilon^2}.\end{equation}
Loosely speaking, equation \eqref{eq:appgeogintro} describes a
physical system where the current point $m(t)$ moves by inertia in
$\mM_N$, but is deflected by a spring of strength $\frac{1}{\epsilon}$
attached to the nearest point $s(t)$ in $\mS$.  Note that the squared
distance $\dd^2_{\mS}$ is semi-concave, and that its restriction to
the finite-dimensional space $\mM_N$ is differentiable at almost every
point.

We now rewrite this systems of equations \eqref{eq:appgeogintro} in
terms of projection on the sets $\mS$ and $\mM_N$. Since the space of
measure-preserving maps $\mS$ is closed but not convex, the orthogonal
projection of $\mS$ exists but is usually not uniquely defined. To
simplify the exposition we will nonetheless associate to any point
$m\in \mM$ one of its projection $P_\mS(m)$, i.e. any point in $\mS$
such that $\nr{P_\mS(m)-m}_\mM = \dd_\mS(m)$. We also denote
$P_{\mM_N}: \mM\to\mM_N$ the orthogonal projection on the linear
subspace $\mM_N\subseteq \mM$. We can rewrite
Eq. \eqref{eq:appgeogintro} in terms of these two projection
operators:
\begin{equation}
\begin{cases}
  \ddot m(t) + \frac{m(t) - P_{\mM_N}\circ P_\mS(m(t))}{{\epsilon^2}} =0, &\hbox{ for } t>0  \, ,  \medskip\\
  (m(0),\dot{m}(0)) \in  \mM_N^2
\end{cases} \label{eq:appgeog2intro}
\end{equation}
From Proposition~\ref{prop:numerics}, the double projection
$P_{\mM_N}\circ P_\mS(m)$ is uniquely defined for almost every $m\in
\mM_N$.  We first prove that the system of equations
\eqref{eq:appgeogintro} can be used to approximate regular solutions
to Euler's equation \eqref{eq:EI}.  Our proof of convergence uses a
modulated energy technic and requires a Lipschitz regularity
assumption on the solution of Euler's equation. It also requires a
technical condition on the computational domain.
\begin{Def}
  An open subset $\Omega$ of $\Rsp^d$ is called prox-regular with
  constant $r_\Omega > 0$ if every point within distance $r_\Omega$
  from $\Omega$ has a unique projection on $\Omega$.
\end{Def}
Note that smooth and semi-convex domains are prox-regular for a
constant $r_\Omega$ smaller than the minimal curvature radius of the
boundary $\partial \Omega$. On the other hand, convex domains are
prox-regular with constant $r_\Omega = +\infty$.

\begin{Thm}\label{Thm:SDreg}
Let $\Omega$ be a connected prox-regular set. Let $v,p$ be a strong
solution of Euler's equations \eqref{eq:EI}, let $\phi$ be the flow map induced by $v$ given by \eqref{eq:EIEtoL} and assume that $v, p, \partial_t v, \partial_t p,
\nabla v$ and $\nabla p$ are Lipschitz on $\Omega$, uniformly on $[0,T] $.
Suppose in addition that there exist a $\Class^1$ curve $m:[0,T]
\to\Rsp$ satisfying the initial conditions
$$ m(0) = P_{\mM_N}(\id),\quad \dot{m}(0) = P_{\mM_N}(v(0,\cdot)),$$
which is twice differentiable and satisfies the second-order equation
\eqref{eq:appgeogintro} for all times in $[0,T]$ up to a (at most)
countable number of exceptions.  Then,
$$ \max_{t\in [0,T]} \nr{\dot m - v(t, \phi(t,\cdot))}^2_\mM \leq C_1 \frac{h_N^2}{\eps^2} + C_2 \eps^2 + C_3 h_N $$
where the constants $C_1$, $C_2$ and $C_3$ only depend on the proximal
constant of the domain, on the $\LL^{\infty}$ norm (in space) of the velocity $v(t,\cdot)$ and on the Lipschitz norms
(in space) of the velocity and its first derivatives
$v(t,\cdot),\nabla v(t,\cdot), \partial_t v(t,\cdot)$ and of the
pressure and its derivatives $p(t,\cdot),\nabla p(t,\cdot), \partial_t
p(t,\cdot)$.
\end{Thm}


The value of $C_1$, $C_2$ and $C_3$ is given more precisely at the end of
Section~\ref{sec:convergencesemidiscret}. Note that the hypothesis on
the solution $m$ to the EDO is here for technical reasons. Removing it
was not of our main concern in this paper since we also give a proof
of convergence of the fully discrete numerical scheme without this
assumption. It is likely that solutions to the EDO
\eqref{eq:appgeogintro} satisfying this hypothesis can be constructed
through di Perna-Lions or Bouchut-Ambrosio theory
\cite{ambrosio2004transport,bouchut2001renormalized,lions1998equations}.

\subsection*{Discretization in space and time}
To obtain a numerical scheme we also need to discretize in time the
Hamiltonian system \eqref{eq:appgeog2intro}. For simplicity of the
analysis, we consider a symplectic Euler scheme.  Let $\tau$ be the
time step, for $m \in \mM_N$ we denote by $P_{\mM_N}P_{\mS}(m) $ a
random element in this set. The solution is the set of points
$M^{n},V^n$ given by:
\begin{equation} \begin{cases} 
(M^0, V^0) \in \mM_N \\
V^{n+1}=V^n -\frac{\tau}{\epsilon^2}  \left( M^{n} -P_{\mM_N}\circ P_{\mS}(M^{n}) \right)  \\
M^{n+1}=M^{n}+\tau V^{n+1} 
\end{cases}    \label{schemeBrenier}
\end{equation}
We also set $t^n = n \tau$. For the numerical scheme of our
approximate geodesic flow we set a more precise theorem.
\begin{Thm}\label{Thm:SDregd}
Let $\Omega$ be a connected prox-regular set, $\epsilon$ and $\tau$ be positive
numbers. Let $v,p$ be a strong solution of \eqref{eq:EI}, let $\phi$ be the flow map induced by $v$ given by \eqref{eq:EIEtoL} and assume that $v, p, \partial_t v, \partial_t p,
\nabla v$ and $\nabla p$ are Lipschitz on $\Omega$, uniformly on $[0,T] $. Let $\left(M^{n},V^n \right)_{n\geq 0}$ be a sequence
generated by \eqref{schemeBrenier} with initial conditions
$$ M^0 = P_{\mM_N}(\id),~~ V^0 = P_{\mM_N}(v(0,\cdot)).$$ Then,
$$ \max_{n\in \N \cap [0,T/\tau]} \nr{V^n - v(t^n,\phi(t^n,\cdot))}_\mM \leq C(h_N \epsilon^{-1}, \tau \epsilon^{-2} ) \left[   \epsilon^2 +h_N + \frac{h^2_N}{ \epsilon^2} + \frac{ \tau}{ \epsilon^2} \right ] $$
where the constant $C$ only depends on upper bounds of $\tau \epsilon^{-2}$  and $h_N \epsilon^{-1}$, on the proximal
constant of the domain, on the $\LL^{\infty}$ norm (in space) of the velocity $v(t,\cdot)$ and on the Lipschitz norms
(in space) of the velocity and its first derivatives
$v(t,\cdot),\nabla v(t,\cdot), \partial_t v(t,\cdot)$ and of the
pressure and its derivatives $p(t,\cdot),\nabla p(t,\cdot), \partial_t
p(t,\cdot)$.

\end{Thm}


In order to use the numerical scheme \eqref{schemeBrenier}, one needs
to be able to compute the double projection operator $P_{\mM_N}\circ
P_{\mS}$ or equivalently the gradient of the squared distance
$\dd_\mS^2$ for (almost every) $m$ in $\mM_N$.  Brenier's polar
factorization problem \cite{brenier1991polar} implies that the squared
distance between a map $m:\Omega\to\Rsp$ and the set $\mS$ of
measure-preserving maps equals the squared Wasserstein distance
\cite{oldandnew} between the restriction of the Lebesgue measure to
$\Omega$, denoted $\Leb$, and its pushforward $m_\# \Leb$ under the
map $m$:
$$ \dd_\mS^2(m) = \min_{s\in \mS} \nr{m - s}^2 =
\Wass_2^2(m_\#\Leb,\Leb). $$ Moreover, since $m$ is piecewise-constant
over the partition $(\omega_i)_{1\leq i\leq N}$, the push-forward
measure $m_\#\Leb$ if finitely supported. Denoting $M_i
\in \Rsp^d$ the constant value of the map $m$ on the subdomain
$\omega_i$ we have,
$$ m_\#\Leb = \sum_{1\leq i\leq N} \Leb(\omega_i)\delta_{M_i} =
\frac{1}{N} \sum_{1\leq i\leq N} \delta_{M_i}.$$ Thus, computing the
projection operator $P_{\mS}$ amounts to the numerical resolution of
an optimal transport problem between the Lebesgue measure on $\Omega$
and a finitely supported measure. Thanks to recent work
\cite{aurenhammer1998minkowski,merigot2011multiscale,de2012blue,levy2014numerical}, this
problem can be solved efficiently in dimensions $d=2,3$. We give more
details in Section~\ref{sec:numerics}.

\begin{Rk}
A scheme involving similar ideas, and in particular the use of optimal
transport to impose incompressibility contraints, has recently been
proposed for CFD simulations in computer graphics \cite{Pixar}. From
the simulations presented in \cite{Pixar}, the scheme seems to behave
better numerically, and it also has the extra advantage of not
depending on a penalization parameter $\eps$. It would therefore be
interesting to extend the convergence analysis presented in
Theorem~\ref{Thm:SDregd} to the scheme presented in \cite{Pixar}. This
might however require new ideas, as our proof techniques rely heavily
on the fact that the space-discretization is hamiltonian, which does
not seem to be the case in \cite{Pixar}.
\end{Rk}

\begin{Rk}
  Our discretization \eqref{eq:appgeogintro} resembles (and is
  inspired by) a space-discretization of Euler's equation
  \eqref{eq:EI} introduced by Brenier in
  \cite{brenier2000derivation}. The domain is also decomposed into
  subdomains $(\omega_i)_{1\leq i\leq N}$, and one considers the set
  $\mS_N\subseteq \mS$, which consists of measure-preserving maps
  $s:\Omega\to\Omega$ that are induced by a permutation of the
  subdomains. Equivalently, one requires that there exists
  $\sigma:\{1,\hdots,N\}\to \{1,\hdots,N\}$ such that $s(\omega_i) =
  \omega_{\sigma(j)}$. The space-discretization considered in
  \cite{brenier2000derivation} leads to an ODE similar to
  \eqref{eq:appgeogintro}, but where the squared distance to $\mS$ is
  replaced by the squared distance to $\mS_N$. This choice of
  discretization imposes strong contraints on the relative size of the
  parameters $\tau$, $h_N$ and $\epsilon$, namely that $h_N =
  \mathrm{O}(\eps^8)$ and $\tau = \mathrm{O}(\eps^4)$. Such
  constraints still exist with the discretization that we consider
  here, but they are milder. In Theorem \ref{Thm:SDregd} the condition $\tau =o(\epsilon^{2})$ is due to the time discretization of \eqref{eq:appgeog2intro} and can be improved using a scheme more accurate on the conservation of the Hamiltonian \eqref{Hamilto}. However even with an exact time discretization of the Hamiltonian, the condition $\tau =o(\epsilon)$ remains mandatory, see section \ref{sec:schemeBrenier}.
\end{Rk}


\subsection*{Acknowledgements}
We would like to thank Yann Brenier for pointing us to
\cite{brenier2000derivation}, and for many interesting discussions at
various stages of this work.

\section{Preliminary discussion on geodesics}\label{sec:basics}
To illustrate the approached geodesic scheme we focus on the very
simple example of $\R$ seen as $\R \times \{0\} \subset \R^2$. The
geodesic is given by the function $\gamma$: $[0,T] \to \R^2$ with
\begin{equation}
\begin{cases}
\gamma(t)=(t,0), \, t \in [0,T], \\
\gamma (0) = (0,0),\\
\dot \gamma (0) = (1,0).\\
\end{cases}   \label{geodesicR}
\end{equation}
We suppose that we make an error of order $(h_0,h_1)$ in the initial conditions.  As in \eqref{eq:appgeogintro}  we consider the solutions of the Hamiltonian system associated to:
\begin{equation}\label{HdansR}
H(z,v)=\frac12 || v ||^2 + \frac{1}{2\epsilon^2} d_{\R \times \{0\}}^2(z).
\end{equation}
That is 
\begin{equation}
\begin{cases}
\ddot z(t)= \frac{1}{\epsilon^2} \left( P_{\R}(z) -z \right) = \frac{1}{2\epsilon^2} \nabla d_{\R \times \{0\}}^2(z),  \, t \in [0,T], \\
 z(0) = (0,h_0),\\
\dot z (0) = (1,h_1).\\
\end{cases}   \label{appgeodesicR}
\end{equation}
where $P_{\R}(z)$ is the orthogonal projection from $\R^2$ onto $\R \times \{0\}$. Notice that we made a mistake of order $h_0$ on the initial position and $h_1$ on the initial velocity. 
In this case the solution is explicit and reads 
\begin{equation}  \label{appgeodesicRexp}
z(t) = \left(t,h_0 \cos\frac{t}{\epsilon}+\epsilon h_1\sin\frac{t}{\epsilon} \right). 
\end{equation}
A convenient way to quantify how far $z$ is from being a geodesic is to use a modulated energy related to the Hamiltonian $H$ and the solution $\gamma$. 
We define $E_{\gamma}$ by
\begin{equation} \label{modr}
E_{\gamma}(z)= \frac12 || \dot z -\dot \gamma ||^2 + \frac{1}{2\epsilon^2} d_{\R \times \{0\}}^2(z).
\end{equation}
Notice that $E_{\gamma}(z)$ is symmetric since for all $t\in [0,T]$, $d_{\R \times \{0\}}^2(z)=0$. 
A direct computation leads to 
\begin{equation} \label{modrest}
E_{\gamma}(z)= \frac{h_0^2}{\epsilon^2}+ h_1^2.
\end{equation}
This estimates shows that the velocity vector field $\dot z$ converges
towards the geodesic velocity vector fields $ \dot \gamma$ as soon as
$h_0$ goes to $0$ quicker then $\epsilon$. Our construction of
approached geodesics for the Euler equation follow this
idea. Estimates \eqref{modrest} suggests that our convergence results
for the incompressible Euler equation in Theorem \ref{Thm:SDreg} is
sharp. A computation of the Hamiltonian \eqref{HdansR} evaluated on the solution of the Euler symplectic scheme, with $h_1=0$ leads to
$$
H(Z^n,V^n) \leq  (1-\frac{\tau^2}{\epsilon^2} )^{n} H(Z^0,V^0). 
$$
It suggests again that the estimation $ \tau = o(\epsilon^2)$ in Theorem \ref{Thm:SDregd} is sharp, even if one can hope for compensation to have in practice a much better convergence. 



\section{Convergence of the approximate geodesics model}\label{sec:convergencesemidiscret}

\subsection{Preliminary lemma} Before proving Theorem~\ref{Thm:SDreg}, we collect a few useful lemmas.


\begin{Lem}[Projection onto the measure preserving maps $\mS$]
  \label{lem:brenier}
    Let $m \in \mM = \LL^2(\Omega,\Rsp^d)$. There exists a convex
    function $\varphi: \Omega\to\Rsp$, which is unique up to an
    additive constant, such that $s$ belongs to $\Pi_\mS(m)$ if and
    only if $m = \nabla \varphi \circ s$ up to a negligible
    set. Moreover, $m-s$ is orthogonal to $\mdiv \circ s$:
  \begin{equation} \label{orthobre}
    \forall v \in \mdiv, \int_{\Omega} \sca{m(x)-s(x)}{v(s(x))} \dd x =0.
  \end{equation}
\end{Lem}
\begin{proof}
    The first part of the statement is Brenier's polar factorization
    theorem \cite{brenier1991polar}, and the uniqueness of $\phi$ follows from
    the connectedness of the domain.  Using a regularization argument
    we deduce the orthogonality relation
\begin{equation*}
\int_{\Omega} m(x)v(s(x))=\int_{\Omega}\nabla \varphi \circ s(x)v(s(x)) =\int_{\Omega}\nabla \varphi v(x) = -\int_{\Omega} \varphi \nabla \cdot v(x) = 0. \qedhere
\end{equation*}
\end{proof}



\begin{Lem}[Projection onto the piecewise constant set $\mM_N$]\label{lem:merigot}
  The projection of a function $g \in L^2(\Omega,\R^d)$ on $\mM_N$ is
  the following piecewise constant function :
  $$ \Pi_{\mM_N}(g) = \sum_{i=1}^N G_i \One_{\omega_i}, \hbox{ with }
  G_i := \frac{1}{\Leb(\omega_i)} \int_{\omega_i} g(x) \dd
  x$$ and where $\One_{\omega_i}$ is the indicator function of
  the subdomain $\omega_i$.
\end{Lem}

\begin{proof}
  It suffices to remark that for any $m\in \mM_N$, $m = \sum_{1\leq
    i\leq N} M_i \One_{\omega_i}$,
  \begin{align}\label{orthomer}
 \sca{g}{m}_\mM &= \int_\Omega \sca{m(x)}{g(x)}\dd x
    = \sum_{1\leq i\leq N} \sca{M_i}{\int_{\omega_i} g(x) \dd x}
    = \sca{m}{\sum_i G_i \One_{\omega_i}}_\mM
    \qedhere  \notag
  \end{align}
\end{proof}

%
%

\begin{Lem}\label{extension}
  Let $\Omega$ be a prox-regular domain of $\Rsp^d$ let $(V,\nr{.})$
  be a normed vector space.  Then, there exists a linear map $L:
  \Class^0(\Omega, V) \to \Class^0(\Rsp^d, V)$ such that for any $f\in
  \Class^0(\Omega,V)$,
  \begin{enumerate}
  \item[(i)] $\restr{L f}_{\Omega} = f$ and $\nr{Lf}_{\LL^{\infty}(\Rsp^d,V)} \leq  \| f \|_{\LL^{\infty}(\Omega,V)}$
  \item[(ii)]  $\Lip(Lf) \leq \frac{6}{r_\Omega} \nr{f}_{\LL^{\infty}(\Omega,V)} + \Lip(f)$.
  \end{enumerate}
\end{Lem}

\begin{proof} Let $r_\Omega$ be the prox-regularity constant of $\Omega$, and let
  $\Omega'$ be a tubular neighborhood of radius $r_\Omega/2$ around
  $\Omega$, i.e. $\Omega' = \{ x\in \Rsp^d \mid \dd(x,\Omega) \leq
  r_\Omega/2 \}.$ Denote $p: \Omega'\to \bar \Omega$ the function which
  maps a point of $\Omega'$ to the closest point in $\bar \Omega$. From
  Theorem 4.8.(8) in \cite{federer1959curvature}, the map $p$ is
  $2$-Lipschitz. We now define the function $Lf$ by
  \begin{equation*}
    Lf(x) = \begin{cases}
      \chi(\nr{x - p(x)})f(p(x)) &\hbox{ if } x \in \Omega'\\
      0  &\hbox{ if not, }
    \end{cases}
  \end{equation*}
  where $\chi(r) = \max(1- 2r/r_\Omega, 0)$. Remark that $\| \chi
  \|_{L^{\infty}}$ is bounded by one, implying that $\| Lf
  \|_{L^{\infty}(\Rsp^d)} \leq \| f \|_{L^{\infty}(\Omega)}$.  For the
  Lipschitz continuity estimates we distinguish three cases.  First,
  if $x,y$ both belong to $\Omega' \times \Omega'$, we have
  \begin{align*}
    \nr{Lf(x)- Lf(y)}&=  \nr{\chi(\nr{x - p(x)})f(p(x)) - \chi(\nr{y - p(y)})f(p(y))}  \\
  &\leq  |\chi(\nr{x - p(x)})|\cdot\nr{f(p(x))-f(p(y))} \\
    & \qquad \qquad+ |\chi(\nr{y - p(y)})-\chi(\nr{x - p(x)})|\cdot \nr{f(p(y))}\\
  &\leq
  \| \chi \|_{\LL^{\infty}(\Rsp)}\Lip(f) \cdot \nr{x-y}
  + \Lip(\chi \circ (\Id-p))\| f \|_{L^{\infty}(\Omega)} \nr{x-y} \\
  &\leq \left(\frac{6}{r_\Omega}\| f \|_{L^{\infty}(\Omega)}+\Lip(f) \right) \nr{x-y}.
  \end{align*}
  If $x$ belongs to $\Omega'$ and $y$ belongs to $\Rsp^d\setminus \Omega'$ one has $Lf(y) = 0$ so that
  \begin{align*}
    \nr{Lf(x)- Lf(y)} &\leq  \nr{\chi(\nr{x - p(x)}) f(p(x))} \\
    &\leq \frac{ 2}{r_\Omega} \| f \|_{L^{\infty}(\Omega)}  \abs{r_\Omega/2 - \nr{x - p(x)}} \\
 &  \leq \frac{ 2}{r_\Omega} | \| f \|_{L^{\infty}(\Omega)}  \nr{y-x}. 
  \end{align*}
  Finally, if $x,y$ are outside of $\Omega'$, $Lf(x) = Lf(y) = 0$ and
  there is nothing to prove.
 \end{proof}

\label{sec:semidiscret}
We are now ready to prove Theorem \ref{Thm:SDreg}. In the following
the dot refer to the time derivative and $\sca{.}{.}$ to the Hilbert
scalar on $\mM$. By abuse of notation we denote by the same variable a
Lipschitz function defined on $\Omega$ and its (also Lipschitz)
extension defined on the whole space $\Rsp^d$ thanks to Lemma
\ref{extension}. The space $\Rsp^d$ is equipped with the Euclidian
norm, and the space of $d\times d$ matrices are equiped with the dual norm.
The Lipschitz constants that we consider are with respect to these two
norms. Finally for a curve $X: t \in [0,T] \mapsto X(t,\cdot)$ we denote
$\LipT(X)= \sup_{t\in[0,T]} \Lip(X(t,\cdot))$.
\subsection{Proof of Theorem \ref{Thm:SDreg}}
Let $v$ be a solution of \eqref{eq:EI} and $m$ a solution of
\eqref{eq:appgeogintro} and for any $t\in[0,T]$, denote $s(t) =
P_{\mS}(m(t))$. In other words, $s(t)$ is an arbitrary choice of a
projection of $m(t)$ on $\mS$.  Equation \eqref{eq:appgeogintro} is the ODE associated to the Hamiltonian
$$H(m,v)=\frac{1}{2}\|v\|^2_{\mM}+\frac{d^2_{\mS}(m)}{2\epsilon^2}.$$
We therefore consider a energy involving this Hamiltonian, modulated
with the exact solution $v$:
\begin{equation}\label{modulatedbrenier}
E_v(t)= \frac{1}{2}\| \dot m(t) - v(t,m(t))\|^2_{\mM}+ \frac{d^2_{\mS}(m)}{2\epsilon^2}.
\end{equation}
We will control $E_v$ using a Gronwall estimate.
\begin{Rk}
Note that we need to use Lemma \ref{extension} to define the modulated
energy $E_v$ since the maps $m(t,\cdot) \in \mM_N$ can send points
outside of $\Omega$ when $\Omega$ is not convex.
\end{Rk}
\subsubsection{Time derivative}
We compute $\frac{d}{dt} E_v(t) $ and modify the expression in order
to identify terms of quadratic order. Since the Hamiltonian $H(\dot
m(t),m(t))$ is preserved, we find
\begin{equation}
\frac{d}{dt} E_v(t) = \underbrace{ - \left< \ddot m(t), v(t,m(t))  \right>}_{I_1}  \underbrace{- \left< \dot m(t) - v(t,m(t)), \partial_t v(t,m(t)) + \left( \dot m(t) \cdot \nabla \right)v(t,m(t))  \right>.}_{I_2} 
\end{equation}
Using the EDO \eqref{eq:appgeogintro}, $I_1$ can be rewritten as
\begin{align*}
{\epsilon^2} I_1 &=   \left< m(t) - P_{\mM_N}(s(t)) , v(t,m(t)) \right> \\
&=  \left< m(t) - s(t) , v(t,m(t)) \right>+  \left< s(t) - P_{\mM_N}(s(t)) , v(t, m(t))\right>\\
 &=  \underbrace{  \left< m(t) - s(t) ,v(t,m(t))- v(t,s(t))\right>}_{{\epsilon^2} I_3},
\end{align*}
where we have used that $s(t) - P_{\mM_N}(s(t))$ is orthogonal to
$\mM_N$ and that $m(t) - s(t)$ is orthogonal to $\mdiv \circ s$, see
Lemmas \ref{lem:merigot} and \ref{lem:brenier}.
To handle the term $I_2$ we define for $X\in \mM$ the two following
operators, often called  material derivatives:
\begin{equation} \label{materialderivative}
\begin{cases}
D_t v (t,X) &= \partial_t v(t,X) + \left( v(t,X) \cdot \nabla \right)v(t,X),\\
D_t p(t, X)&= \partial_t p(t,X)+ \left< v(t,X) ,\nabla p (t,X))  \right>.
 \end{cases}
 \end{equation}
Remark that Euler's equation \eqref{eq:EI} implies that $D_t v
(t,s(t))=- \nabla p (t,s(t))$. This leads to 
\begin{align*}
I_2 &=  - \left< \dot m(t) - v(t,m(t)), \partial_t v(t,m(t)) + \left( v(t,m(t)) \cdot \nabla \right)v(t,m(t))  \right> \\
& \underbrace{- \left< \dot m(t) - v(t,m(t)),  \left( \dot m(t) - v(t,m(t)) \cdot \nabla \right)v(t,m(t))  \right>}_{I_4} \\
&= I_4   \underbrace{- \left< \dot m(t) - v(t,m(t)), D_t v(t,m(t)) - D_t v(t,s(t))  \right> }_{I_5}  + \underbrace{\left< \dot m(t) - v(t,m(t)) ,\nabla p (t,s(t)) \right>}_{I_6}    \\
\end{align*}
We rewrite $I_6$  as
\begin{align*}
I_6&= \underbrace{- \left< \dot m(t)-v(t,m(t)),\nabla p (t,m(t)) - \nabla p (t,s(t)) \right>}_{I_7} +  \left< \dot m(t)-v(t,m(t)),\nabla p (t,m(t)) \right> \\
&=  I_7 + \frac{d}{dt} \underbrace{ \int_{\Omega} p(t,m(t,x)))dx}_{-J(t)}  -  \int_{\Omega} \partial_t p(t,m(t,x))-   \left< v(t,m(t,x)) ,\nabla p (t,m(t,x))  \right>dx \\
&=-\frac{d}{dt} J(t)+I_7 - \underbrace {\int_{\Omega} D_t p(t,m(t,x) )dx}_{I_8}.
\end{align*}
%
%

\begin{Rk}
  The quantity $I_5+I_7$ control the fact that the extension of
  $(v,p)$ constructed by Lemma~\ref{extension} is not a solution of
  the Euler equation on $\Rsp^d$ (in particular, $I_5+I_7$ vanishes if
  $(v,p)$ solves Euler's equation on $\Rsp^d$).
\end{Rk}

\subsubsection{Estimates}
Many of the integrals $I_i$ can be easily bounded using the energy
$E_v$ and Cauchy-Schwarz and Young's inequalities. First, 
\begin{align}
\nonumber I_3&\leq \left| \frac{ \left< m(t) - s(t) ,v(t,m(t))- v(t,s(t))\right>}{\epsilon^2} \right| \\
 \label{estimationI3} &\leq  \Lip( v(t) )  \frac{\| m(t) - s(t) \|_{\mM}^2 }{\epsilon^2}   \leq \LipT( v ) E_v(t). 
\end{align}
Furthermore 
\begin{equation}\label{estimationI4}
I_4 \leq \sup_{x\in \Rsp^d} || \nabla v(t,x) ||  \| \dot m(t) - v(t,m(t))\|_{\mM}^2 \leq   \LipT (v) E_v(t),
\end{equation}
Where $C$ depends only on the dimension $d$. 
To estimate $I_5$ and later $I_8$ we first remark that $D_tv$ and $D_tp$ are Lipschitz operators with 
\begin{equation} \label{bornelipv} 
  \begin{aligned}
  \LipT(D_tv)&\leq \LipT( \partial_t v )+ \LipT(v)  \|  \nabla v \|_{L^{\infty}} + \LipT(\nabla v)\| v \|_{L^{\infty}}\\
  &\leq \LipT( \partial_t v )+ \LipT(v)  \LipT (v) + \LipT(\nabla v)\| v \|_{L^{\infty}}
  \end{aligned}
\end{equation}
\begin{equation}
  \begin{aligned}
    \LipT(p_tv)&\leq \LipT( \partial_t p )+ \LipT(v) \|  \nabla p \|_{L^{\infty}} + \LipT(\nabla p)   \| v \|_{L^{\infty}} \\
\label{bornelipp}&\leq \LipT( \partial_t p )+ \LipT(v) \LipT(p)  + \LipT(\nabla p)\| v \|_{L^{\infty}}.
  \end{aligned}
\end{equation}
For $I_5$ we obtain --- using 
$\dd_\mS(m(t)) = \nr{m(t) - s(t)}_\mM \leq \epsilon \sqrt{E_v(t)} $ and $ \| \dot m(t) - v(t,m(t))  \|_{\mM}\leq \sqrt{E_v(t)}  $ to get from the second to the third line ---,
\begin{align}
\nonumber I_5 &\leq  \left|\left< \dot m(t) - v(t,m(t)), D_t v(t,m(t)) - D_t v(t,s(t))  \right> \right| \\
\nonumber&\leq  \LipT(D_tv) \| \dot m(t) - v(t,m(t))  \|_{\mM} \| m(t) - s(t)  \|_{\mM} \\
\label{estimationI5}&\leq \epsilon \LipT(D_tv) E_v(t)
\end{align}
The quantity $I_7$ can be bounded using the same arguments,
\begin{align}
\nonumber I_7 &\leq \left| \left< \dot m(t)-v(t,m(t)),\nabla p (t,m(t)) - \nabla p (t,s(t)) \right> \right| \\
\nonumber&\leq  \LipT(\nabla p) \| \dot m(t) - v(t,m(t))  \|_{\mM} \| m(t) - s(t)  \|_{\mM} \\
\label{estimationI7}&\leq  \epsilon  \LipT(\nabla p) E_v(t).
\end{align}
Finally to estimate $I_8$ and $J$ we can assume that $\int_{\Omega}
p(t,x)\dd x = 0$ since the pressure is defined up to a constant. Using that $s(t)$ is measure-preserving, this gives 
\begin{align*}
\int_{\Omega} D_t p(t, s(t,x))dx &= \int_{\Omega} \partial_t p(t,s(t,x)) + \left< v(t,s(t,x)) ,\nabla p (t,s(t,x)))  \right>dx\\
&= \int_{\Omega} \partial_t p(t,x))dx + \int_{\Omega} \left< v(t,x) ,\nabla p (t,x)) \right>dx  =0,
\end{align*}
Therefore, using Young's inequality,
\begin{align}
\nonumber  I_8 &\leq  \left| \int_{\Omega} D_t p(t, m(t,x))dx - \int_{\Omega} D_t p(t, s(t,x))dx \right| \leq \LipT(D_t p) \| m(t)-s(t) \|_{L^1(\Omega)} \\
\nonumber &\leq \frac12 \frac{|| m(t)-s(t)||^2_{L^2{(\Omega)}}}{2\epsilon^2}+C  \LipT(D_t p) \epsilon^2\\
\label{estimationI8} &\leq  \frac12 E_v(t)+ \cst(\Omega) \LipT(D_t p) \epsilon^2,
\end{align}
where in this estimates and in the following estimates $\cst(\Omega)$
is a constant depending only on $\Leb(\Omega)$. Similarly
\begin{align}
\nonumber |J(t)|&\leq \left| \int_{\Omega} p(t,m(t,x))) -p(t,s(t,x))) dx \right| \leq \LipT(p) || m(t)-s(t)||_{L^1(\Omega)} \\
\label{estimationJ} &\leq \frac{1}{2} E_v(t)+ \cst(\Omega)  \LipT(p)\epsilon^2.
\end{align}
Remark also that 
\begin{equation}\label{estimationJ0}
\abs{J(0)}\leq  \LipT(p) h_N.
\end{equation}
\begin{Rk}
The two last estimates show that we can add $\frac{\dd}{\dd t}J$ into
the Gronwall argument. It is a general fact that the derivative of a
controlled quantity can be added. This is a classical way of
controlling the term of order one in the energy.
\end{Rk}
\subsection{Gronwall argument}
Collecting estimates \eqref{estimationI3}, \eqref{estimationI4}, \eqref{estimationI5}, \eqref{estimationI7}, \eqref{estimationI8} and \eqref{estimationJ} we get 
\begin{align*}
&\frac{d}{dt} \left( E_v(t) + J(t) \right) \leq I_3+I_4+I_5+I_7+I_8 + J(t) - J(t) \\
&\leq   \left[ 2\LipT( v ) + \epsilon \LipT(D_tv) + \epsilon  \LipT(\nabla p)+ 1   \right] \left(E_v(t)+J(t)\right)\\
& + \cst(\Omega) \left( \LipT(D_t p) + \LipT(p) \right) \epsilon^2
\end{align*}
Setting 
\begin{equation*}
\begin{cases}
\tilde C_1&= 2\LipT( v ) + \epsilon \LipT(D_tv) + \epsilon  \LipT(\nabla p)+ 1  , \\
\tilde C_2&=  \left( \LipT(D_t p) + \LipT(p) \right) ,\\
\end{cases}
\end{equation*}
we obtain 
$$
\frac{d}{dt} \left( E_v(t) + J(t) \right) \leq \cst(\Omega)  \tilde C_1 (E_v(t)+J(t))+ \cst(\Omega) \tilde C_2 \epsilon^2.
$$
We deduce that for any $t\in [0,T]$:
$$
E_v(t) \leq \left(\left(E_v(0) + J(0) \right)+ \cst(\Omega) \tilde C_2T\epsilon^2 \right)e^{\tilde C_1T}  -J(t) 
$$
Using one more time \eqref{estimationJ} we obtain
$$
E_v(t) \leq 2 \left(E_v(0) + \LipT(p)h_N + \cst(\Omega) \tilde C_2T\epsilon^2 \right)e^{\tilde C_1T} +
\cst(\Omega)\LipT(p) \epsilon^2.
$$
Finally using that $$E_v(0) =  \frac{1}{2}\|P_{\mM}(v_0) - v_0\|^2_{\mM} + \frac{d^2_{\mS}(\Id)}{2\epsilon^2} \leq  \frac{h^2_N}{2} + \frac{h^2_N}{2\epsilon^2}  $$ 
and 
\begin{align}
\nonumber \nr{\dot m(t) - v(t, \phi(t))}^2_\mM &\leq 2 \nr{\dot m(t) - v(t, m(t))}^2_\mM + \nr{v(t, m(t)) - v(t, \phi(t))}^2_\mM \\
\nonumber & \leq  2 E_v(t) + 2(\LipT(v))^2 \nr{m(t) - \phi(t)}^2_\mM \\ \nonumber  &\leq 2 E_v(t)  + 2(\LipT(v))^2 \dd^2_\mS(m(t))  \\ 
\label{areutiliser}& \leq 2(1+(\LipT(v))^2\epsilon^2)E_v(t).
\end{align}

we conclude 
\begin{align}
\nonumber \nr{\dot m(t) - v(t, \phi(t))}^2_\mM &\leq  (2+(\LipT(v))^2\epsilon^2)E_v(t) \\
\label{closelook} & \leq 2(1+(\LipT(v))^2\epsilon^2)  \left[2 \left( \frac{h^2_N}{2} + \frac{h^2_N}{2\epsilon^2} + \LipT(p)h_N 
\nonumber +  \cst(\Omega)\tilde C_2 T\epsilon^2 \right)e^{\tilde C_1T} \right. \\ &+  \left. \cst(\Omega)  \LipT(p) \epsilon^2\right] \\
\label{cetruc}  &\leq C_1 \frac{h^2_N}{\epsilon^2} + C_2 \epsilon^2 +C_3 h_N
\end{align}
where 
\begin{equation*}
\begin{cases}
C_1 &= 2(1+(\LipT(v))^2\epsilon^2) e^{\tilde C_1T}   \\
C_2 &= \cst(\Omega)(1+(\LipT(v))^2 )\left(   \tilde C_2 T e^{\tilde C_1T} +  \LipT(p)   \right) \\
C_3 &=  \cst(\Omega)(1+(\LipT(v))^2\epsilon^2)\LipT(p)  \\
\end{cases}
\end{equation*}
where we used that $\epsilon$ and $h_N$ are smaller than $ \cst(\Omega)$. Observe that the RHS of \eqref{cetruc} goes to zero as  $\frac{h_N}{\epsilon}$ and $\epsilon$ goes to zero. It finishes the proof of Theorem \ref{Thm:SDreg}.
In order to track down the regularity assumptions, we give the value of $\tilde C_1$, $\tilde C_2$ in term of the data:

\begin{align*}
\tilde C_1&= 1+2\LipT( v )    + \epsilon  \LipT(\nabla p)\\
&+ \epsilon \left( \LipT( \partial_t v )+ (\LipT(v))^2 + \LipT(\nabla v)\| v \|_{L^{\infty}} \right),\\
%
\tilde C_2&=  \left( \LipT(D_t p) + \LipT(p) \right) \\
&=  \LipT(p) + \LipT( \partial_t p )+ \LipT(v) \LipT(p)  + \LipT(\nabla p)\| v \|_{L^{\infty}}. 
\end{align*}

\begin{Rk}\label{rk:gen}
A close look to the explicit value of $\tilde C_1$, $\tilde C_2$ and estimation \eqref{closelook}, together with a diagonal argument shows that our scheme approximate solutions less regular than supposed in Theorem \ref{Thm:SDreg}. For example we can set the following theorem:
Let $v,p$ be a solution of Euler's equation \eqref{eq:EI}. Suppose that $v$ is Lipschitz in space. Suppose although that there exists $(v_k,p_k)_{k\in \N}$ a sequence of regular (in the sense of Theorem \ref{Thm:SDreg}) solutions  of \eqref{eq:EI}
such that $v_k(0,\cdot) \longrightarrow v(0,\cdot) $ in $\mM$ and $\Lip_T(v_k) \longrightarrow \Lip_T(v) $. Then there exists $h_N(k)$ and $\epsilon(k)$, polynomials in the data, such that $ \nr{\dot m(t)[v_k(0),\epsilon(k)] - v\left(t, m(t)[v_k(0),\epsilon(k)]\right)}^2_\mM $ goes to zero as $k$ goes to infinity.  
\end{Rk}

\section{Convergence of the Euler symplectic numerical scheme \label{sec:schemeBrenier}}

In this section we prove Theorem \ref{Thm:SDregd}. The proof follows the one given in \ref{sec:semidiscret} for Theorem \ref{Thm:SDreg} with some additional terms. It combined two Gronwall estimates. The first one is a continuos Gronwall argument on the segment $[n \tau ,(n+1)\tau]$, the second one is a discrete Gronwall argument. For both steps we use the modulated energy.

For a solution of \eqref{schemeBrenier} and $\theta \in [0,\tau]$ we denote


\begin{equation} \begin{cases} \label{scheme}
V^{n+\theta}&=V^n - \theta   \frac{ M^n -P_{\mM}\circ P_{\mS}(M^n) }{\epsilon^2}  \\
M^{n+\theta}&=M^n+ \theta  V^{n+1} ,
\end{cases}   
\end{equation}
the linear interpolation between $(M^n,V^n) $ and $(M^{n+1},V^{n+1} )$. 

\subsection{The modulated energy}
The Hamiltonian at a step $n$ is $$H^n=H(M^n,V^n)=\frac{1}{2}\|V^n\|^2_{\mM}+\frac{d^2_{\mS}(M^n)}{2\epsilon^2}.$$
The modulated energy at time  $n\tau$ is 
\begin{equation}\label{modulatedbrenierd}
E^n= \frac{1}{2}\left\| V^n - v\left(n\tau,M^n\right)\right\|^2_{\mM}+ \frac{d^2_{\mS}(M^n)}{2\epsilon^2}.
\end{equation}
For $\theta \in [0,\tau]$ we consider
\begin{equation}\label{modulatedbreniertheta}
\begin{cases}
H^{n+\theta}&= \frac{1}{2}\left\| V^{n+\theta} \right\|^2_{\mM}+ \frac{\dd^2_{\mS}(M^{n+\theta})}{2\epsilon^2},\\
E^{n+\theta}&= \frac{1}{2}\left\| V^{n+\theta} - v\left(n\tau+\theta,M^{n+\theta}\right)\right\|^2_{\mM}+ \frac{\dd^2_{\mS}(M^{n+\theta})}{2\epsilon^2}.
\end{cases}
\end{equation} 
Remark that 
\begin{equation}\label{modulatedbrenieravecH}
E^{n+\theta}=H^{n+\theta}- \left< V^{n+\theta}, v\left(n\tau+\theta,M^{n+\theta}\right)\right> +  \frac{1}{2}\left\| v\left(n\tau+\theta,M^{n+\theta}\right)\right\|^2_{\mM}
\end{equation}
We start with a lemma quantifying the conservation of the Hamiltonian. 
\begin{Lem}[Conservation of the Hamiltonian]\label{conshami}
For $\theta \in [0,\tau]$ and $n\in \N \cap [0,T/\tau]$ there holds
\begin{equation}\label{hamiltonestimates2}
\left( 1- \frac{\tau ^2}{\epsilon^2} \right)H^{n+1} \leq   H^n,
\end{equation} 

\begin{equation}\label{hamiltonestimates2Bornes}
H^{n} \leq  e^{{T\tau}\epsilon^{-2}}\left(\frac12 \|V^0 \|^2_{\mM}+ \frac{h^2_N}{2\epsilon^2}\right),
\end{equation}
and 
\begin{equation}\label{hamiltonestimates}
H^{n+\theta} \leq H^n + C({\tau}\epsilon^{-2},h_N\epsilon^{-1}) \frac{\tau^2}{\epsilon^2}, 
\end{equation}
where $C({\tau}\epsilon^{-2},h_N\epsilon^{-1})$ depends only on $\|V^0 \|_{\mM}^2$, $T$ and upper bounds of ${\tau}\epsilon^{-2}$ and $h_N\epsilon^{-1}$.

\end{Lem}
\begin{proof}
The proof is based on the  $\frac{1}{\epsilon^2}$-semiconcavity of $\frac{d^2_{\mS}}{2}$, see Proposition \ref{prop:numerics} for details.
On the one hand the  $\frac{1}{\epsilon^2}$-semiconcavity of $\frac{d^2_{\mS}}{2}$ reads
\begin{equation*}
\frac{d^2_{\mS}(M^{n+\theta})} {2\epsilon^2} \leq  \frac{d^2_{\mS}(M^{n})}{2\epsilon^2} + \theta \left<  V^{n+1} ,   \frac{ M^n -P_{\mM}\circ P_{\mS}(M^n) }{\epsilon^2}  \right> + \frac{\theta^2}{2\epsilon^2} \|  V^{n+1} \|^2_{\mM},
\end{equation*}
where we used that $  \left[ M^n -P_{\mM}\circ P_{\mS}(M^n) \right] \in \nabla \dd^2_{\mS}(M^{n})  $ and \eqref{scheme}.  On the other hand, \eqref{scheme} again, leads
$$
\frac{ \| V^{n+\theta}\|^2_{\mM}}{2} = \frac{\| V^{n}\|^2_{\mM}}{2} - {\theta} \left<  V^{n} ,   \frac{ M^n -P_{\mM}\circ P_{\mS}(M^n) }{\epsilon^2}  \right> + {\theta^2} \left\|  \frac{ M^n -P_{\mM}\circ P_{\mS}(M^n) }{\epsilon^2}\right \|^2_{\mM}
$$
Summing both equations and using \eqref{scheme} gives
\begin{equation}\label{pdd}
H^{n+\theta} \leq H^n+ \frac{\theta(\tau-\theta)}{\epsilon^2} \frac{ \left \| M^n -P_{\mM}\circ P_{\mS}(M^n) \right \|^2_{\mM} }{\epsilon^2} +  \frac{\theta ^2}{\epsilon^2} \frac{\| V^{n+1} \|^2_{\mM}}{2}. 
\end{equation}
Applied with $\theta=\tau$, it proves \eqref{hamiltonestimates2}. The inequality \eqref{hamiltonestimates2Bornes} is a direct consequence of \eqref{hamiltonestimates2}. To obtain \eqref{hamiltonestimates} remark that by definition of the projection $P_{\mS}$ 
\begin{align}
\nonumber \left\|  M^n - P_{\mM}\circ P_{\mS}(M^n) \right\|_{\mM} & \leq 2 \left\|  M^n -s^n\right\|_{\mM} + 2 \left\|s^n - P_{\mM}\circ P_{\mS}(M^n) \right\|_{\mM}\\
\label{gradphipluspetitquephi} & \leq 2 \left\|  M^n - s^n \right\|_{\mM}  = 2 d_{\mS}(M^{n}).
\end{align}

Therefore \eqref{pdd} rewrites
\begin{align*}
H^{n+\theta}
&\leq H^n + \frac{\theta^2}{\epsilon^2} H^{n+1} + \frac{8\theta(\tau-\theta)}{\epsilon^2}  H^n.\\
\end{align*}
Combined with \eqref{hamiltonestimates2Bornes}, it proves \eqref{hamiltonestimates} and finishes the proof of Lemma \ref{conshami}.
\end{proof}

We deduce from \eqref{modulatedbrenieravecH} and \eqref{hamiltonestimates} that for any $\theta \in [0,\tau]$ and ${n\in \N \cap [0,T/\tau]}$
\begin{equation}\label{modulatedbrenieravecH2}
E^{n+\theta} \leq E^n + \int_0^1 d^{n+\theta} d\theta + C({\tau}\epsilon^{-2},h_N\epsilon^{-1}) \frac{\tau^2}{\epsilon^2} 
\end{equation}
where $$d^{n+\theta}=\frac{d}{d\theta}\left[- \left< V^{n+\theta}, v\left(n\tau+\theta,M^{n+\theta}\right)\right> +  \frac{1}{2}\left\| v\left(n\tau+\theta,M^{n+\theta}\right)\right\|^2_{\mM} \right].$$
We compute $d^n$ using \eqref{scheme}
and the notations $v^{n+\theta}_{p}= v(n\tau+\theta,M^p)$, $s^{n+\theta}= P_{\mS}(M^{n+\theta})$ and  $v^{n+\theta}_{s^{n+\theta}}= v(n\tau+\theta,s^{n+\theta})$.

\begin{align*}
d^{n+\theta} &= - \left< \frac{d}{d\theta} V^{n+\theta}, v^{n+\theta}_{n+\theta} \right> 
- \left< V^{n+\theta}, \frac{d}{d\theta} v^{n+\theta}_{n+\theta})\right>   +  \left< v^{n+\theta}_{n+\theta}, \frac{d}{d\theta}  v^{n+\theta}_{n+\theta}  \right> \\
&=  \underbrace{ \epsilon^{-2} \left<M^n -P_{\mM}\circ P_{\mS}(M^n)   ,v^{n+\theta}_{n+\theta} \right>}_{I_1}   \underbrace{- \left< V^{n+\theta} - v^{n+\theta}_{n+\theta} , \partial_t v^{n+\theta}_{n+\theta}+ \frac{d}{d\theta} M^{n+\theta} \cdot \nabla v ^{n+\theta}_{n+\theta}  \right>}_{I_2} 
\end{align*}
The term $I_1$ rewrites
\begin{align*}
{\epsilon^2} I_1 &=  \left<  M^n -P_{\mM}\circ P_{\mS}(M^n) , v^{n+\theta}_{n+\theta}  \right>   \\
&=    \left<  M^n -s^n , v^{n+\theta}_{n+\theta}  \right> +    \left< s^n -P_{\mM} \circ P_{\mS}(M^n) , v^{n+\theta}_{n+\theta} \right> \\ 
&=  \underbrace{  \left<  M^n -s^n , v^{n+\theta}_{n+\theta}  - v^{n+\theta}_{s^n}\right>}_{{\epsilon^2} I_3} 
\end{align*}

Here we had to control the fact that, due to the double projection, the norm of the acceleration $ \left\| M^n -P_{\mM}\circ P_{\mS}(M^n) \right\|^2_{\mM}$ is not equal to $d^2_{\mS}(M^{n})$. We used the orthogonality property of the double projection to control this problem. On the one hand $ s^n -P_{\mM} \circ P_{\mS}(M^n)$ is orthogonal to $\mM$ since it is a linear subspace. On the other hand $M^n -s^n$ is orthogonal to the tangent space of $\mS$ at $s^n$, see Lemma \ref{lem:brenier}.
 \\

To handle $I_2$ we used the material derivatives defined by \eqref{materialderivative}, 
\begin{align*}
I_2 &=  - \left< V^{n+\theta} - v^{n+\theta}_{{n+\theta}} ,  \partial_t v^{n+\theta}_{n+\theta}+ \frac{d}{d\theta} M^{n+\theta} \cdot \nabla v ^{n+\theta}_{n+\theta}  \right>  \\
I_2 &=  -  \left<V^{n+\theta} - v^{n+\theta}_{{n+\theta}} ,  \partial_t v^{n+\theta}_{n+\theta}+  v ^{n+\theta}_{n+\theta} \cdot \nabla v ^{n+\theta}_{n+\theta}  \right> \\
& \underbrace{- \left<V^{n+\theta} - v^{n+\theta}_{{n+\theta}},   \left( \frac{d}{d\theta} M^{n+\theta} -  v ^{n+\theta}_{n+\theta}\right) \cdot \nabla  v ^{n+\theta}_{n+\theta}  \right>}_{I_4} \\
&= I_4  \underbrace{-\left< \frac{d}{d\theta} M^{n+\theta} -  v ^{n+\theta}_{n+\theta} , D_t v^{n+\theta}_{n+\theta} - D_t s^{n+\theta}_{n+\theta}  \right>}_{I_5} + \underbrace{\left< \frac{d}{d\theta} M^{n+\theta} -  v ^{n+\theta}_{n+\theta} ,\nabla p^{n+\theta}_{s^{n+\theta}}  \right>}_{I_6} \\
\end{align*}We rewrite $I_6$:
\begin{align*}
I_6 &=  \underbrace{ \left< \frac{d}{d\theta} M^{n+\theta} -  v ^{n+\theta}_{n+\theta} ,\nabla p^{n+\theta}_{s^{n+\theta}}  - \nabla p^{n+\theta}_{{n+\theta}} \right>}_{I_7} + \left<  \frac{d}{d\theta}M^{n+\theta} -  v ^{n+\theta}_{n+\theta} ,\nabla  p^{n+\theta}_{n+\theta}   \right> \\
&= I_7 +  \left<  \frac{d}{d\theta}M^{n+\theta}  ,\nabla  p^{n+\theta}_{n+\theta}   \right>-  \left<   v ^{n+\theta}_{n+\theta} ,\nabla  p^{n+\theta}_{n+\theta}   \right> \\
&=  I_7 + \frac{d}{d\theta} \underbrace{ \int_{\Omega} p^{n+\theta}_{n+\theta} dx}_{-J(\theta)}  -   \int_{\Omega} \partial_t p^{n+\theta}_{n+\theta}-  \left< v^{n+\theta}_{n+\theta} ,\nabla  p^{n+\theta}_{n+\theta}  \right> dx \\
&=-\frac{d}{d\theta } J(\theta) -  \underbrace {\int_{\Omega} D_t p^{n+\theta}_{n+\theta} dx}_{I_8},
\end{align*}
\subsection{Gronwall estimates on $[n\tau, (n+1)\tau ]$}
Using \ref{scheme} we obtain for $I_3$:
\begin{align}
\nonumber I_3&=\epsilon^{-2} \left<  M^n -s^n ,  v^{n+\theta}_{n+\theta}  - v^{n+\theta}_{s^n} \right>  \\
\nonumber &\leq  \LipT( v )  \frac{\| M^n -s^n \|_{\mM} \| M^{n+\theta} -s^n\|_{\mM} }{\epsilon^2} \\
\nonumber &\leq  \LipT( v )  \frac{\| M^n -s^n \|_{\mM} \| M^{n+\theta} -M^n\|_{\mM} }{\epsilon^2} +    \LipT( v )  \frac{\| M^n -s^n \|_{\mM} \| M^{n} -s^n\|_{\mM} }{\epsilon^2} \\
\nonumber &\leq \LipT( v ) \left(   \frac{\| M^n -s^n \|^2_{\mM} }{\epsilon^2} + \theta \epsilon^{-1} \frac{\| M^n - s^n \|_{\mM} }{\epsilon}\|V^{n+1}\|_{\mM}       \right) \\
\nonumber &\leq     \LipT( v )  (1+  2\theta  \epsilon^{-1}) E^n  +  2  \LipT( v )   \theta \epsilon^{-1}H^{n}\\
\label{estimationI3discrete}  &\leq  C   \LipT( v )  (1+  2\theta  \epsilon^{-1}) E^n +  2   \LipT( v )  C({\tau}\epsilon^{-2},h_N\epsilon^{-1})  \theta \epsilon^{-1}.
\end{align} 
We used \eqref{hamiltonestimates2Bornes} to obtain the last line. 
Since $\frac{d}{d\theta} M^{n+\theta}=V^{n+1}$, $I_4$ rewrites
\begin{align}
\nonumber I_4&=  -   \left<V^{n+\theta} - v^{n+\theta}_{{n+\theta}},   \left( V^{n+1} -  v ^{n+\theta}_{n+\theta}\right) \cdot \nabla  v ^{n+\theta}_{n+\theta}  \right> \\
\nonumber &=   -   \left<V^{n+\theta} - v^{n+\theta}_{{n+\theta}},   \left( V^{n+\theta} -  v ^{n+\theta}_{n+\theta}\right) \cdot \nabla  v ^{n+\theta}_{n+\theta}  \right> \\
\nonumber& -   \left<V^{n+\theta} - v^{n+\theta}_{{n+\theta}},   \left( V^{n+1} -  V^{n+\theta}\right) \cdot \nabla  v ^{n+\theta}_{n+\theta}  \right>\\
\nonumber  &\leq  || \nabla v(n+\theta) ||_{L^{\infty}(\Omega)}   \left\|V^{n+\theta} - v^{n+\theta}_{{n+\theta}} \right\|^2_{\mM} \\
\nonumber& - (\tau -\theta)\epsilon^{-2}   \left<V^{n+\theta} - v^{n+\theta}_{{n+\theta}},   \left(  M^n -P_{\mM}\circ P_{\mS}(M^n)\right) \cdot \nabla  v ^{n+\theta}_{n+\theta}  \right> \\
\nonumber&\leq \LipT (v) E^{n+\theta} -   (\tau-\theta)\epsilon^{-2} \left<V^{n+\theta} - v^{n+\theta}_{{n+\theta}},   \left(  M^n -s^n \right) \cdot \nabla  v ^{n+\theta}_{n+\theta}  \right> \\
\nonumber&\leq \LipT (v) E^{n+\theta} +  (\tau-\theta)\epsilon^{-1} \left\|V^{n+\theta} - v^{n+\theta}_{{n+\theta}}\right\|_{\mM}  \frac{  \left\|  M^n -s^n \right\|_{\mM} }{\epsilon} \\
\nonumber &\leq   \LipT (v) \left( (1+\frac12(\tau-\theta)\epsilon^{-1}  ) E^{n+\theta}  +\frac12(\tau-\theta)\epsilon^{-1}  E^{n}  \right) \\
\label{estimationI4discrete} &\leq  \LipT (v) (1+\frac12(\tau-\theta)\epsilon^{-1}  ) E^{n+\theta} + \LipT (v) \frac12 (\tau-\theta)\epsilon^{-1}   E^{n}
\end{align}
We used that $ \left<V^{n+\theta} - v^{n+\theta}_{{n+\theta}},   \left(  s^n -P_{\mM}\circ P_{\mS}(M^n)\right) \cdot \nabla  v ^{n+\theta}_{n+\theta}  \right>=0$ since $s^n -P_{\mM}\circ P_{\mS}(M^n)$ is orthogonal to $\mM_N$ and  the quantity $\nabla  v ^{n+\theta}_{n+\theta}$ is a symmetric operator from $\mM_N$ to $\mM_N$. At the antepenultimate line we use Young's inequality. The estimates of $I_5$ and $I_7$ are similar to the semi discrete case.

\begin{align}
\nonumber I_5 & \leq  \left|  \left< V^{n+1} -  v ^{n+\theta}_{n+\theta} , D_t v^{n+\theta}_{n+\theta} - D_t s^{n+\theta}_{n+\theta}  \right> \right| \\
\nonumber&\leq  \LipT(D_tv) \left\| V^{n+1} -  v ^{n+\theta}_{n+\theta} \right\|_{\mM} \left\| M^{n+\theta} - s^{n+\theta}  \right\|_{\mM} \\
\nonumber&\leq  \LipT(D_tv) \left[\left\| V^{n+\theta} -  v ^{n+\theta}_{n+\theta} \right\|_{\mM} \left\| M^{n+\theta} - s^{n+\theta}  \right\|_{\mM} +  \left\| V^{n+1} -V^{n+\theta}  \right\|_{\mM} \left\| M^{n+\theta} - s^{n+\theta}  \right\|_{\mM} \right]\\
\nonumber&\leq   \epsilon \LipT(D_tv)E^{n+\theta} + \LipT(D_tv) (\tau-\theta) \frac{\left\|  M^n - P_{\mM}\circ P_{\mS}(M^n)  \right\|_{\mM}}{\epsilon} \frac{ \left\| M^{n+\theta} - s^{n+\theta}  \right\|_{\mM}}{\epsilon} \\
\label{estimationI5discrete}&\leq \left( \epsilon \LipT(D_tv) +(\tau-\theta)  \LipT(D_tv) \right)  E^{n+\theta} + (\tau-\theta)  \LipT(D_tv) E^n. 
\end{align}
We used Young's inequality and \ref{gradphipluspetitquephi}. 
The quantity $I_7$ is of the same kind. 
\begin{align}
\nonumber I_7 &\leq \left| \left< V^{n+1} -  v ^{n+\theta}_{n+\theta} ,\nabla p^{n+\theta}_{s^{n+\theta}}  - \nabla p^{n+\theta}_{{n+\theta}} \right> \right| \\
\label{estimationI7discrete}&\leq  \left( \epsilon \LipT(\nabla p) +(\tau-\theta)  \LipT(\nabla p) \right)  E^{n+\theta} + (\tau-\theta)  \LipT(\nabla p) E^n 
\end{align}

To estimate $J$ and $I_8$ recall that $\int_{\Omega} D_t p(t, s^n(t,x))dx =0$ and we set $\int_{\Omega} p(t,x)dx=0$. 
\begin{align}
\nonumber  I_8  &\leq \LipT( D_t p ) || M^{n+\theta}-s^{n+\theta}||_{L^1(\Omega)} \leq  \LipT( D_t p ) \left( \frac{|| M^{n+\theta}-s^{n+\theta}||_{\mM}}{2\epsilon^2}+ C\epsilon^2\right)\\
\label{estimationI8discrete} &\leq  \LipT( D_t p ) \frac12 E^{n+\theta}+ C \LipT( D_t p ) \epsilon^2
\end{align}
Similarly 
\begin{align}
\nonumber |J(n\tau +\theta) | &\leq \left| \int_{\Omega} p^{n+\theta}_{n+\theta} - p^{n+\theta}_{s^{n+\theta}}  dx \right|  \leq \LipT(p) || M^{n+\theta}-s^{n+\theta}||_{L^1(\Omega)} \\
\label{estimationJdiscrete}  &\leq   \frac12 E^{n+\theta}+C\epsilon^2
\end{align}
\subsection{Gronwall argument on $[n \tau, (n+1)\tau]$}
From now and for clarity we do not track the constants anymore, $C$ will be a constant depending only on $T$, $\Omega$, $\LipT( v )$, $\LipT( \nabla p )$, $\LipT( D_t v )$ and $\LipT( D_t p )$. The constant $C$ can change between estimates. 
Collecting estimates \eqref{estimationI3discrete}, \eqref{estimationI4discrete}, \eqref{estimationI5discrete}, \eqref{estimationI7discrete}, \eqref{estimationI8discrete} and \eqref{estimationJdiscrete} and intergreting $\theta$ from $0$ to $\tau $ we get 
\begin{align}
J^{n+\theta} + \int_0^{\theta} d^{n+s}ds & \leq J^{n} + C \tau (1+  \tau \epsilon^{-1}) E^n +   C({\tau}\epsilon^{-2},h_N\epsilon^{-1})  \tau^2 \epsilon^{-1} \\
\nonumber & + \int_0^{\theta} C (1+\frac12(\tau-s)\epsilon^{-1}  ) E^{n+s}ds +  C \tau^2 \epsilon^{-1}   E^{n} \\
\nonumber &+  2 \int_0^{\theta} C\left( \epsilon +(\tau-s)  \right)  E^{n+s}ds + C \tau^2  E^n \\
\nonumber & +\int_0^{\theta}\frac12 E^{n+s}ds+ C\tau \epsilon^2 + \int_0^{\theta} J^{n+s} -  J^{n+s} ds \\
\label{Gronlocal}&\leq  J^{n} + C \tau (1+  2 \tau \epsilon^{-1}) E^n +  2C\tau \epsilon^2 + C({\tau}\epsilon^{-2},h_N\epsilon^{-1})   \tau^2 \epsilon^{-1}  \\
\nonumber & + C\int_0^{\theta}  (2+2(\tau-s)\epsilon^{-1}  ) \left( E^{n+s} +J^{n+s} \right) ds.  \\
\end{align}
Remark that we only kept the first order terms using $\epsilon \leq C$ thus  $ \tau^2 \leq C \tau^2 \epsilon^{-1}$ and $(\tau-s) \leq C(\tau-s) \epsilon^{-1}$. Plugging \eqref{Gronlocal} into \eqref{modulatedbrenieravecH2} we obtain 
\begin{align*}
E^{n+\theta}+J^{n+\theta}  & \leq  E^n + J^{n} + C({\tau}\epsilon^{-2},h_N\epsilon^{-1}) \tau^2 \epsilon^{-2} \\
&+ C \tau (1+  \tau \epsilon^{-1}) E^n +    C({\tau}\epsilon^{-2},h_N\epsilon^{-1})\tau^2 \epsilon^{-1}  + C\tau \epsilon^2  \\
& + C \int_0^{\theta}  (2+2(\tau-s)\epsilon^{-1}  ) \left( E^{n+s} +J^{n+s} \right) ds. 
\end{align*}
The Gronwall Lemma on $[0,\tau]$ implies 
\begin{align*}
E^{n+\theta}+J^{n+\theta}  & \leq \left[  E^n + J^{n}+ C({\tau}\epsilon^{-2},h_N\epsilon^{-1}) \tau^2\epsilon^{-2}  \right. \\
  & \left. + C \tau (1+ \tau \epsilon^{-1}) E^n +   C({\tau}\epsilon^{-2},h_N\epsilon^{-1})  \tau^2 \epsilon^{-1} + C\tau \epsilon^2 \right] e^{C\tau(1 +\tau \epsilon^{-1}) }.
\end{align*}
and in particular 
\begin{align}
\nonumber E^{n+1}+J^{n+1}  & \leq \left[ \left(1+ C \tau (1+ \tau \epsilon^{-1}) \right) \left( E^n + J^{n} \right) \right. \\
\label{Grondiscret}  &+ \left.  C({\tau}\epsilon^{-2},h_N\epsilon^{-1}) \tau^2\epsilon^{-2}+ C({\tau}\epsilon^{-2},h_N\epsilon^{-1}) \tau^2 \epsilon^{-1} + C\tau \epsilon^2 \right] e^{C\tau(1 +\tau \epsilon^{-1}) } .
\end{align}
\subsection{Discrete Gronwall step}
From \eqref{Grondiscret} and the descrete Gronwall inequality we deduce, for any ${n\in \N \cap [0,T/\tau]}$:
\begin{multline}
E^{n}+J^{n}  \leq \\ C  \left[  E^0 + J^{0}  + TC({\tau}\epsilon^{-2},h_N\epsilon^{-1})  \left( \tau \epsilon^{-1}+\tau \epsilon^{-2}\right) +T \epsilon^2 \right]  e^{ T(1 + \tau \epsilon^{-1} ) } e^{CT}e^{CT \tau \epsilon^{-1} }.
\end{multline}
Using once again \eqref{estimationJdiscrete} leads
$$
E^{n}  \leq  C\left[  E^0 + J^{0}   +  C({\tau}\epsilon^{-2},h_N\epsilon^{-1})  \left( \tau \epsilon^{-1}+\tau \epsilon^{-2}\right) +  \epsilon^2 \right] e^{ \tau\epsilon^{-2} }e^{\tau \epsilon^{-1}  } + C\epsilon^2.
$$
Including the initial error and rearranging the terms yields 
$$
E^{n}  \leq  C(h_N \epsilon^{-1}, \tau \epsilon^{-2} ) \left[   \epsilon^2 +h_N + \frac{h^2_N}{ \epsilon^2} + \frac{ \tau}{ \epsilon} + \frac{ \tau}{ \epsilon^2}\right ]. 
$$
Using \eqref{areutiliser} we conclude
\begin{multline}
 \max_{n\in \N \cap [0,T/\tau]} \nr{V^n - v(t^n,\phi(t^n,\cdot))}_\mM \\ \leq  C  \max_{n\in \N \cap [0,T/\tau]} E^{n}  \leq C(h_N \epsilon^{-1}, \tau \epsilon^{-2} ) \left[   \epsilon^2 +h_N + \frac{h^2_N}{ \epsilon^2} + \frac{ \tau}{ \epsilon}+ \frac{ \tau}{ \epsilon^2} \right ]. 
\end{multline}
It finishes the proof of Theorem \ref{Thm:SDregd}.
\begin{Rk}
A close look to the constant leads to a similar result as the one given in Remark \ref{rk:gen}: namely the convergence of the numerical scheme towards less regular solutions of the Euler's equations. 
\end{Rk}
\begin{Rk}
The condition $\tau=o(\epsilon^2)$ is linked to the estimate on the Hamiltonian \eqref{hamiltonestimates} in \eqref{hamiltonestimates2} and precisely arises in Lemma \ref{conshami}. Another time discretization, with a better estimate at this stage, would improve this condition. However the bounds in Lemma \eqref{conshami} seems very pessimistic. Experimentally,  the Hamiltonian seems very-well preserved and therefore the convergence criteria is more likely to be $\tau=o(\epsilon)$.
\end{Rk}


\section{Numerical implementation and experiments}
\label{sec:numerics}

\subsection{Numerical implementation}
We discuss here the implementation of the numerical scheme
\eqref{schemeBrenier} and in particular the computation of the double
projection $P_{\mM_N}\circ P_{\mS}(m)$ for a piecewise constant
function $m\in \mM_N$. Using Brenier's polar factorisation theorem,
the projection of $m$ on $\mS$ amounts to the resolution of an optimal
transport problem between $\Leb$ and the finitely supported measure
$m_\#\Leb$. Such optimal transport problems can be solved numerically
using the notion of Laguerre diagram from computational geometry.
\begin{Def}[Laguerre diagram]
  Let $M = (M_1,\hdots,M_N) \in (\Rsp^d)^N$ and let $\psi_1,\hdots,\psi_N \in
  \Rsp$.  The Laguerre diagram is a decomposition of $\Rsp^d$ into
  convex polyhedra defined by
  $$ \Lag_i(M, \psi) = \left\{ x \in \Rsp^d \mid \forall j \in \{1,\hdots,N\},~~
 \nr{x - M_i}^2 + \psi_i \leq \nr{x - M_i}^2 + \psi_j \} \right\}. $$
\end{Def}
In the following proposition, we denote $\Pi_\mS(m) = \{ s \in \mS
\mid \nr{m - s} = \dd_\mS(m) \}.$

\begin{Prop}\label{prop:numerics}
    Let $m \in \mM_N \setminus \mD_N$ and define $M_i = m(\omega_i)
    \in \Rsp^d$. There exist scalars $(\psi_i)_{1\leq i\leq N}$, which
    are unique up to an additive constant, such that
    \begin{equation}
    \label{eq:DMA}
    \forall i\in\{1,\hdots,N\},\qquad\Leb(\Lag_i(M,\psi)) = \frac{1}{N} \Leb(\Omega)
    \end{equation}
    We denote $L_i := \Lag_i(M,\psi)$.Then, a function $s \in \mS$ is
    a projection of $m$ on $\mS$ if and only if it maps the subdomain
    $\omega_i$ to the Laguerre cell $L_i$ up to a negligible set, that
    is:
    \begin{equation}
       \Pi_\mS(m) = \{ s \in \mS \mid \forall i\in\{1,\hdots,N\},~~\Leb(s(\omega_i) \Delta L_i) = 0 \} 
    \end{equation}
    where $\Delta$ denotes the symmetric difference.
    Moreover, $\dd^2_\mS(m)$ is differentiable at $m$ and, setting $B_i = \frac{1}{\Leb(L_i)}\int_{L_i} x \dd x$,
    \begin{equation}
    \begin{aligned}
      \dd^2_\mS(m) &= \sum_{1\leq i\leq N} \int_{L_i} \nr{x - M_i}^2 \dd x\\
      \nabla \dd^2_\mS(m) &= 2 (m - P_{\mM_N} \circ P_\mS(m)) \hbox{ with }
      P_{\mM_N} \circ P_\mS(m) = \sum_{1\leq i \leq N} B_i \One_{L_i}.
    \end{aligned}
    \end{equation}
\end{Prop}

\begin{proof}
  The existence of a vector $(\psi_i)_{1\leq i\leq N}$ satisfying
  Equation~\eqref{eq:DMA} follows from optimal transport theory (see
  Section 5 in \cite{aurenhammer1998minkowski} for a short proof), and
  its uniqueness follows from the connectedness of the domain
  $\Omega$. In addition, the map $T:\Omega \to \{ M_1,\hdots, M_N \}$
  defined by $T(L_i) = M_i$ (up to a negligible set) is the gradient
  of a convex function and therefore a quadratic optimal transport
  between $\Leb$ and the measure $\frac{\Leb(\Omega)}{N} \sum_i
  \delta_{M_i}$. By Brenier's polar factorization theorem, summarized in Lemma~\ref{lem:brenier},
  \begin{align*}
    s \in \Pi_{\mS}(m) \Longleftrightarrow m = T \circ s \hbox{ a.e.} 
    &\Longleftrightarrow \forall i \in \{1,\hdots,N\},~~ \Leb(\omega_i \Delta (T\circ s)^{-1}(\{M_i\})) = 0 \\
    &\Longleftrightarrow \forall i \in \{1,\hdots,N\},~~ \Leb(s(\omega_i) \Delta L_i) = 0,
  \end{align*}
  where the last equality holds because $s$ is measure preserving. To
  prove the statement on the differentiability of $\dd_\mS^2$, we
  first note that the function $\dd_\mS^2$ is $1$-semi-concave, since
  $$ D(m) := \nr{m}^2 - \dd^2_\mS(m) = \nr{m}^2 - \min_{s \in \mS}
  \nr{m - s}^2 = \max_{s \in \mS} 2\sca{m}{s} - \nr{s}^2$$ is
  convex. The subdifferential of $D$ at $m$ is given by $\partial D(m)
  = \{ P_{\mM_N}(s) \mid s \in \Pi_{\mS}(m) \},$ so that $D$ (and
  hence $\dd_{\mS}^2$) is differentiable at $m$ if and only if
  $P_{\mM_N}(\Pi_\mS(m))$ is a singleton. Now, note from
  Lemma~\ref{lem:merigot} that for $s \in \Pi_{\mS}(m)$
  $$ P_{\mM_N}(s) = \sum_{1\leq i\leq N} \bary(s(\omega_i))
  \One_{\omega_i} = \sum_{1\leq i\leq N} \bary(L_i)
  \One_{\omega_i}. $$ This shows that $P_{\mM_N}(\Pi_\mS(m))$ is a
  singleton, and therefore establishes the differentiability of $\dd_\mS^2$ at
  $m$, together with the desired formula for the gradient.
\end{proof}

The difficulty to implement the numerical scheme \eqref{schemeBrenier}
is the resolution of the discrete optimal transport problem
\eqref{eq:DMA}, a non-linear system of equations which must be solved
at every iteration. We resort to the damped Newton's algorithm
presented in \cite{kitagawa2016newton} (see also
\cite{mirebeau2015discretization}) and more precisely on its
implementation in the open-source PyMongeAmpere
library\footnote{\url{https://github.com/mrgt/PyMongeAmpere}}.

\subsubsection{Construction of the tessellation of the domain}
\label{subsec:optimomega}
The fixed tessellation $(\omega_i)_{1\leq i\leq N}$ of the domain
$\Omega$ is a collection of Laguerre cells that are computed through a
simple fixed-point algorithm similar to the one presented in
\cite{de2012blue}. We start from a random sampling $(C^0_i)_{1\leq
  i\leq N}$ of $\Omega$. At a given step $k\geq 0$, we compute
$(\psi_i)_{1\leq i\leq N} \in \Rsp^N$ such that
$$ \forall i \in \{1,\hdots,N\},~ \Leb(\Lag_i(C,\psi)) = \frac{1}{N}
\Leb(\Omega), $$ and we then update the new position of the centers
$(C^{k+1}_i)$ by setting $C^{k+1}_i := \bary(\Lag_i(C^k,\psi)).$ After
a few iterations, a fixed-point is reached and we set $\omega_i :=
\Lag_i(C^k,\psi)$.

\subsubsection{Iterations}
To implement the symplectic Euler scheme for \eqref{eq:appgeog2intro},
we start with $M^0_i := \bary(\omega_i)$ and $V^0_i :=
v_0(M^0_i)$. Then, at every iteration $k\geq 0$, we use Algorithm~1 in
\cite{kitagawa2016newton} to compute a solution $(\psi_{i}^k)_{1\leq
  i\leq N} \in \Rsp^N$ to Equation~\eqref{eq:DMA} with $M = M^k$, i.e.
such that
    $$ \forall i\in\{1,\hdots,N\},~~ \Leb(\Lag_i(M^k,\psi^k)) =
\frac{1}{N} \Leb(\Omega). $$ Finally, we update the positions
$(M^{k+1}_i)_{1\leq i\leq N}$ and the speeds $(V^{k+1}_i)_{1\leq i\leq
  N}$ by setting
    \begin{equation} \label{eq:upd}
      \left\{
    \begin{aligned}
      &V^{k+1}_i = V^k_i +\frac{\tau}{\eps^2}(\bary(\Lag_i(M^k,\psi^k) - M^k_i)\\
      &M^{k+1}_i = M^k_i + \tau V^{k+1}_i
    \end{aligned}
    \right.
    \end{equation}    

\subsection{Beltrami flow in the square}
\begin{figure}
  \includegraphics[width=.32\textwidth]{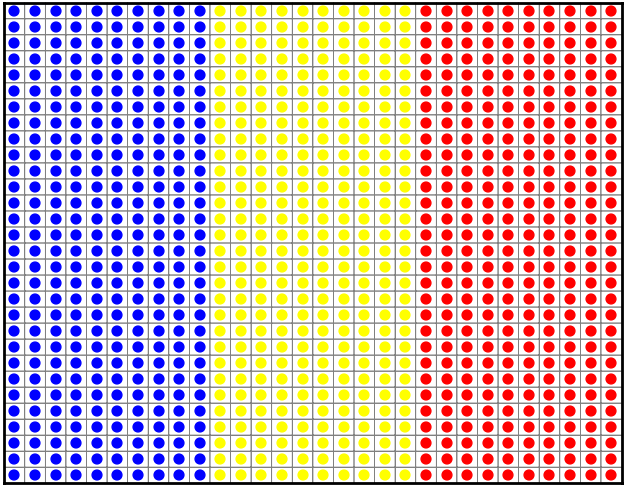}
  \hfill
  \includegraphics[width=.32\textwidth]{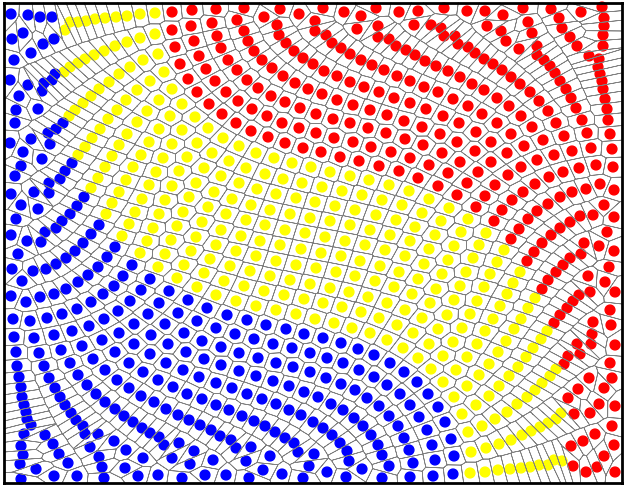} \hfill
  \includegraphics[width=.32\textwidth]{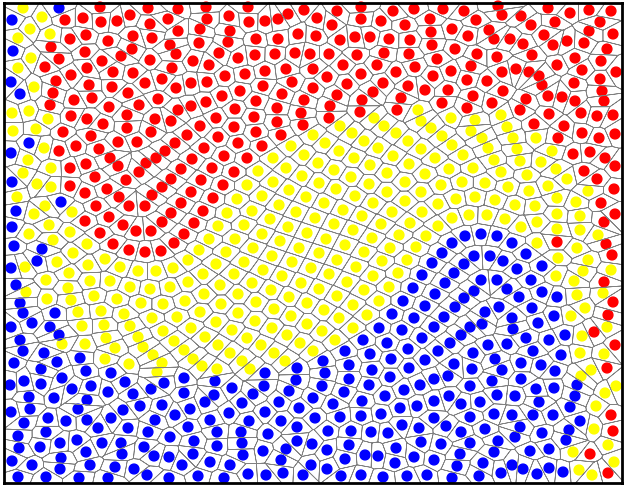} \\
    \includegraphics[width=.32\textwidth]{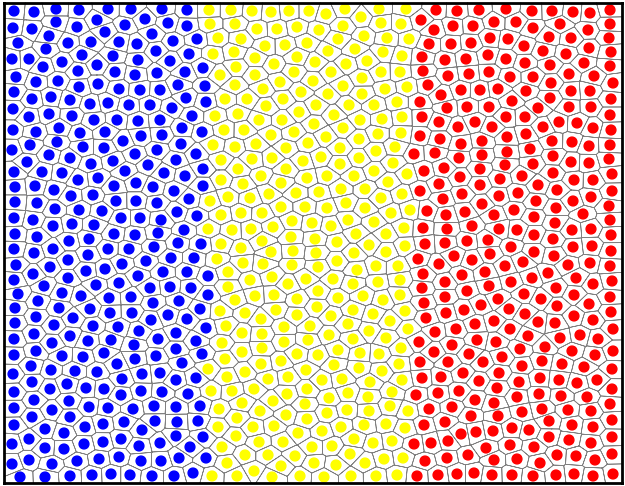} \hfill
  \includegraphics[width=.32\textwidth]{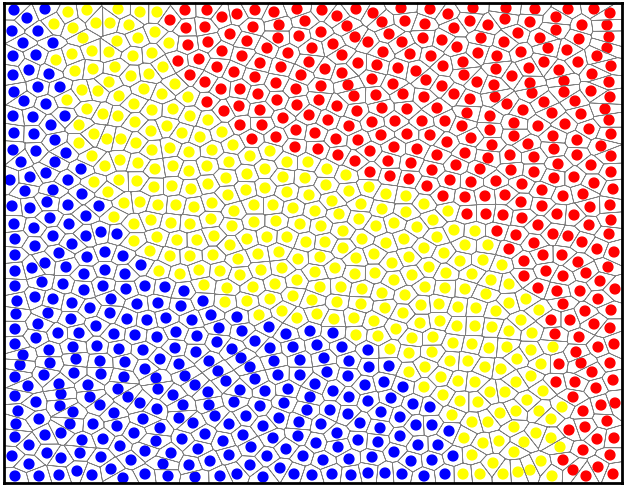} \hfill
  \includegraphics[width=.32\textwidth]{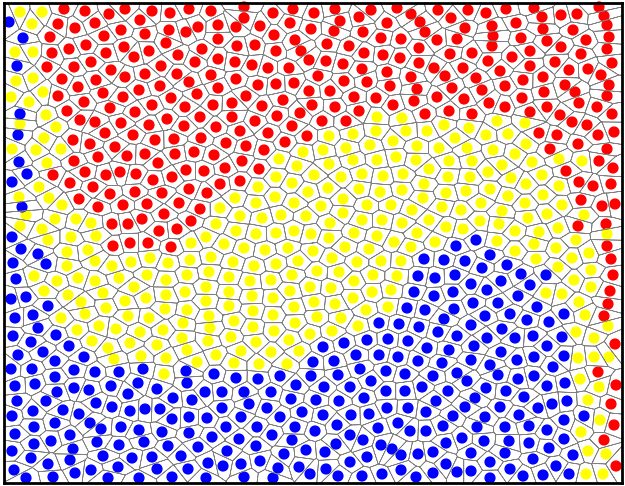} \\
  \caption{(Top row) Beltrami flow in the square, with $N=900$
    particles, $\tau=1/50$ and $\eps = .1$. The particles are colored
    depending on their initial position in the square. From left to
    right, we display the Laguerre cells and their barycenters at
    timesteps $k=0, 24$ and $49$. The partition $(\omega_i)_{1\leq
      i\leq N}$ is induced by a regular grid. (Bottom row) Same
    experiment, but where the partition $(\omega_i)_{1\leq i\leq N}$
    is optimized using the algorithm described in
    \S\ref{subsec:optimomega}. \label{fig:laguerre}}
\end{figure}

Our first testcase is constructed from a stationary solution to
Euler's equation in $2$D. On the unit square $\Omega = [-\frac{1}{2},\frac{1}{2}]^2$, we
consider the Beltrami flow constructed from the time-independent
pressure and speed:
\begin{equation*}
  \left\{
\begin{aligned}
p_0(x_1,x_2) &= \frac{1}{2}(\sin(\pi x_1)^2 + \sin(\pi x_2)^2)\\
v_0(x_1,x_2) &= (-\cos(\pi x_1) \sin(\pi x_2), \sin(\pi x_1) \cos(\pi x_2)) 
\end{aligned}
\right.
\end{equation*}
In Figure~\ref{fig:laguerre}, we display the computed numerical
solution using a low number of particles ($N=900$) in order to show
the shape of the Laguerre cells associated to the solution.

\subsection{Kelvin-Helmoltz instability}
\begin{figure}
  \includegraphics[width=.31\textwidth]{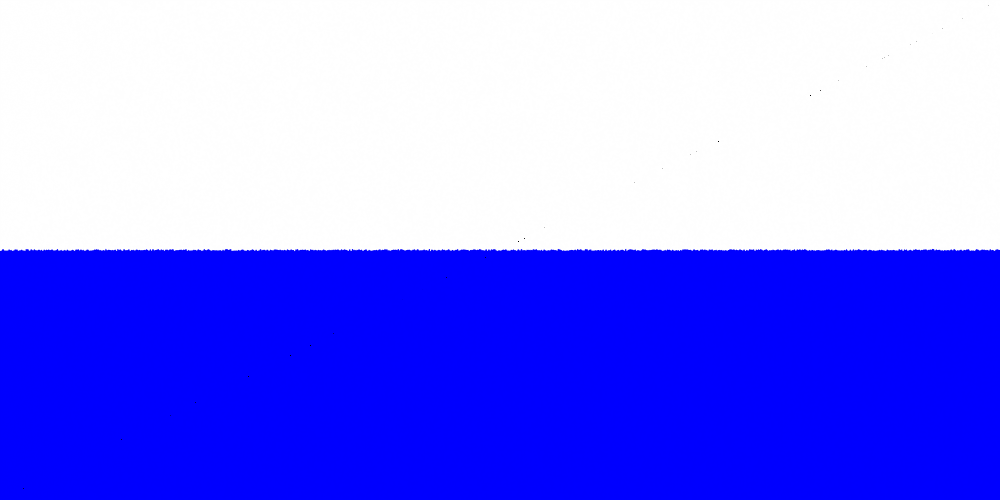} \hspace{.05cm}
  \includegraphics[width=.31\textwidth]{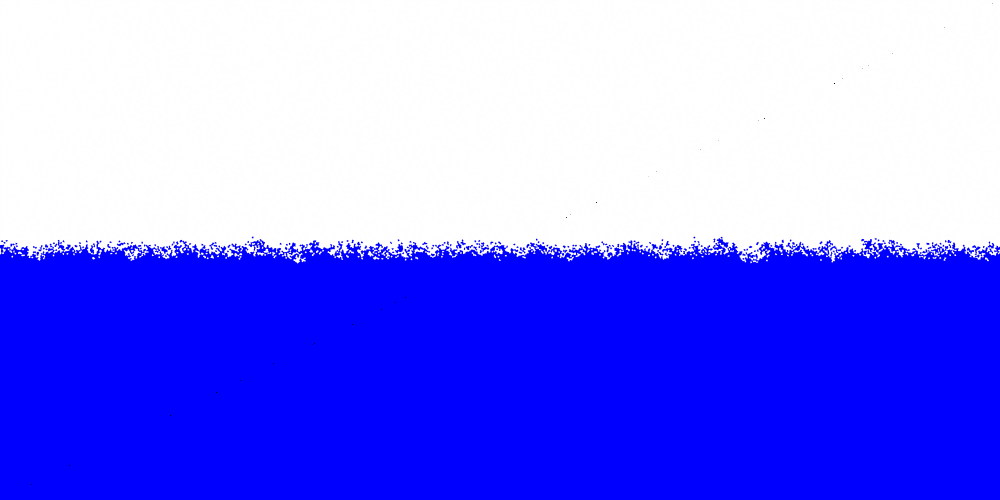} \hspace{.05cm}
  \includegraphics[width=.31\textwidth]{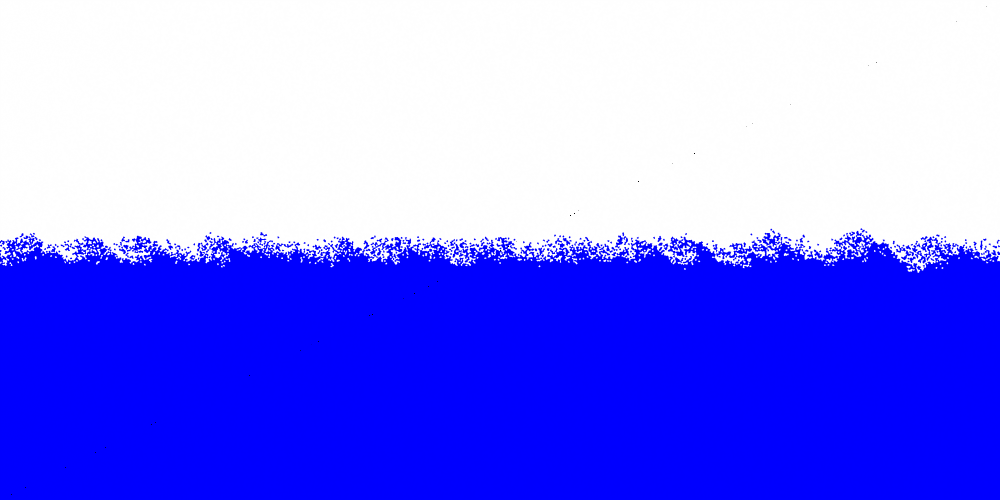} \\
  \includegraphics[width=.31\textwidth]{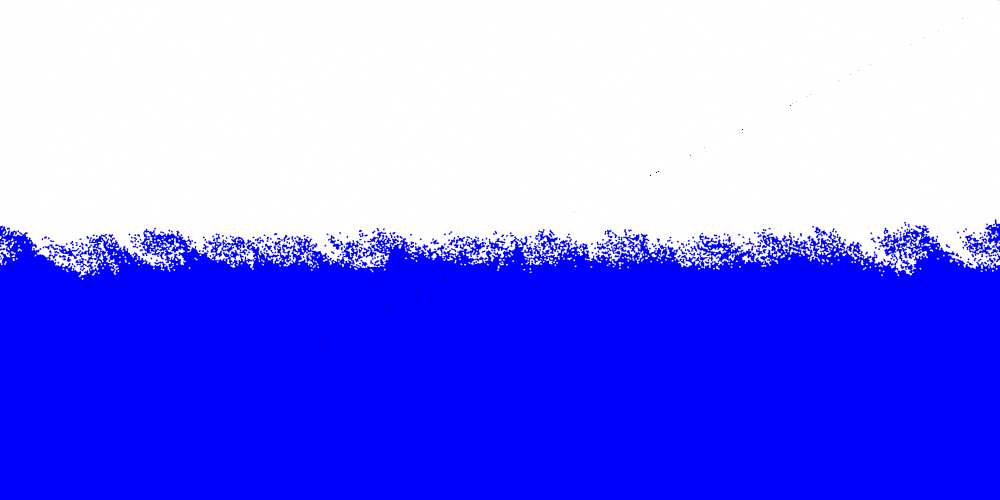} \hspace{.05cm}
  \includegraphics[width=.31\textwidth]{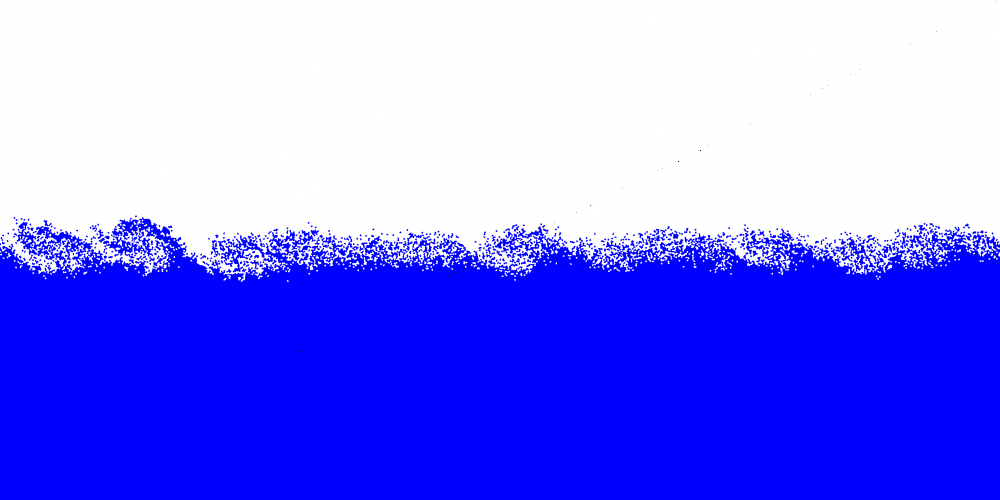} \hspace{.05cm}
  \includegraphics[width=.31\textwidth]{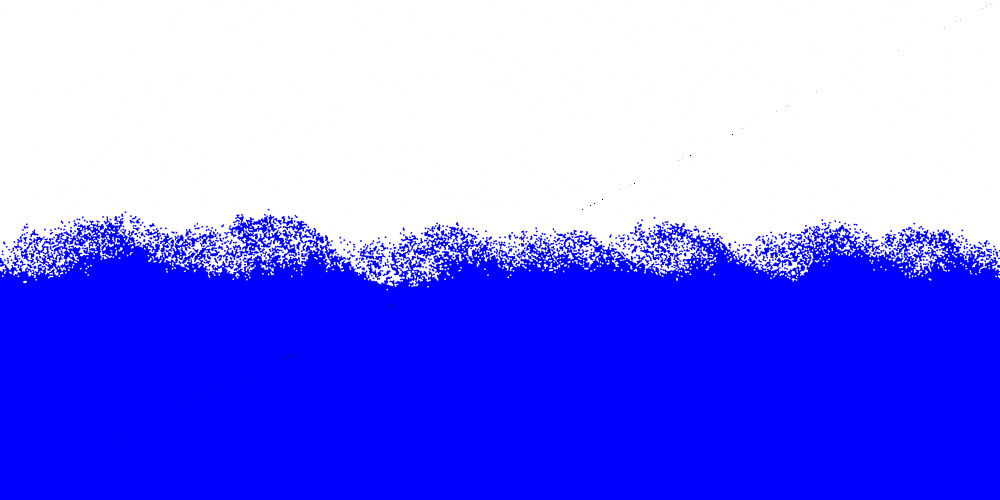} \\
  \includegraphics[width=.31\textwidth]{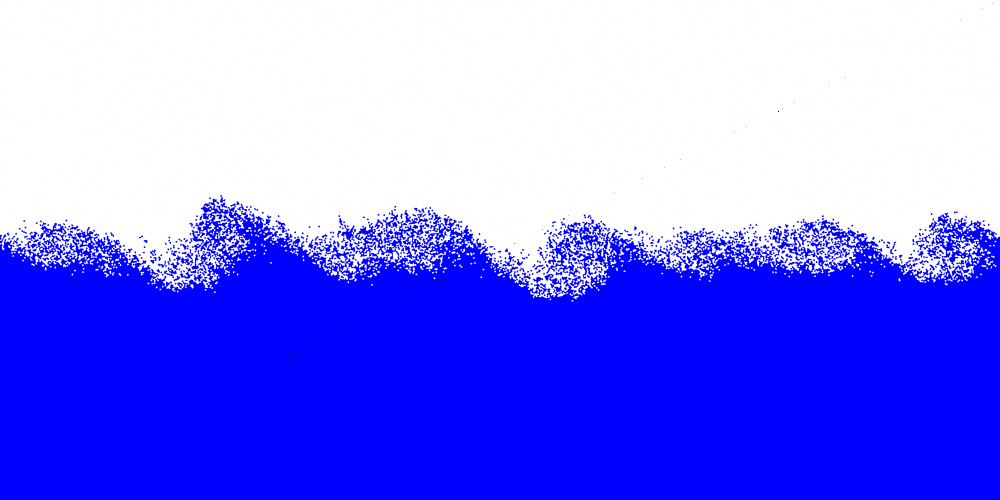} \hspace{.05cm} 
  \includegraphics[width=.31\textwidth]{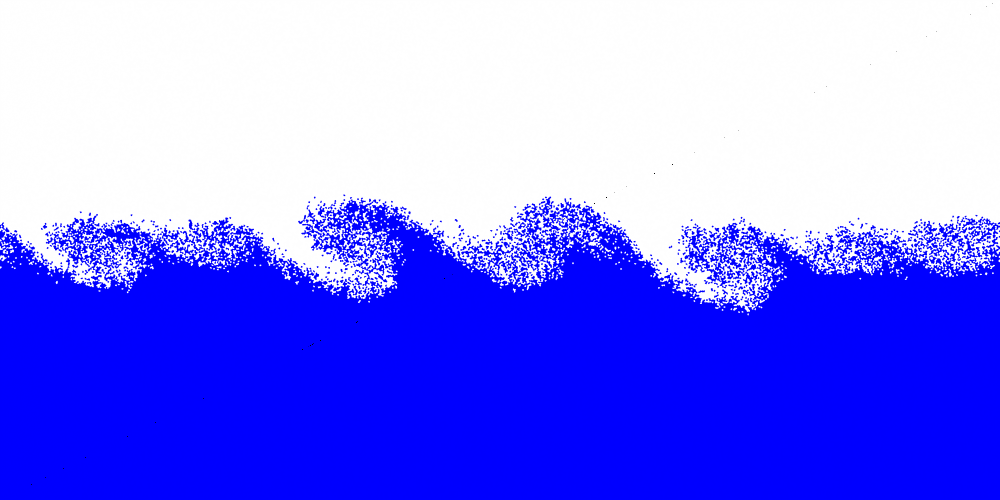} \hspace{.05cm} 
  \includegraphics[width=.31\textwidth]{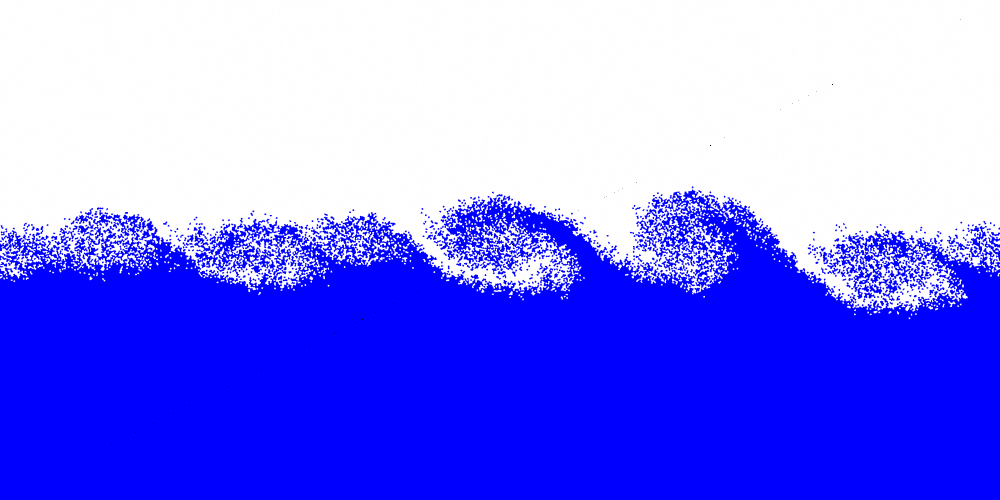} \hspace{.05cm}\\ 
  \includegraphics[width=.31\textwidth]{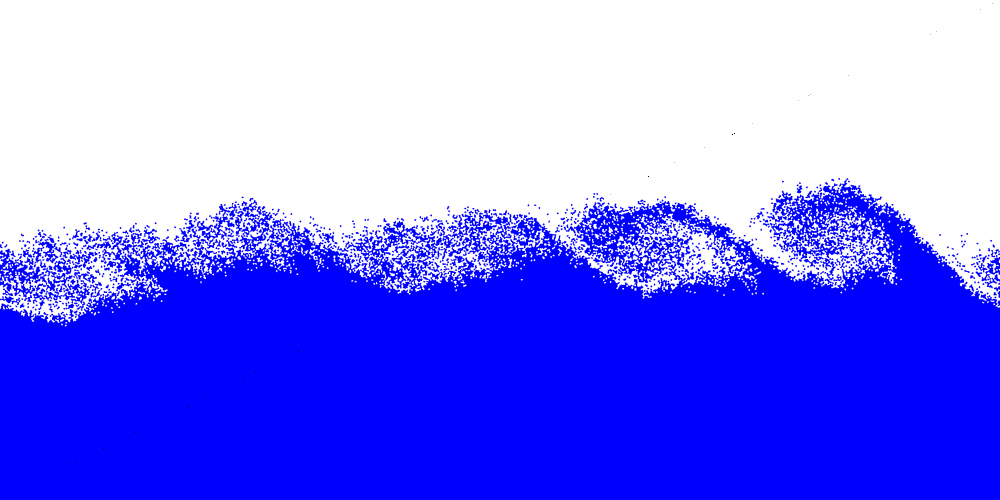} \hspace{.05cm} 
  \includegraphics[width=.31\textwidth]{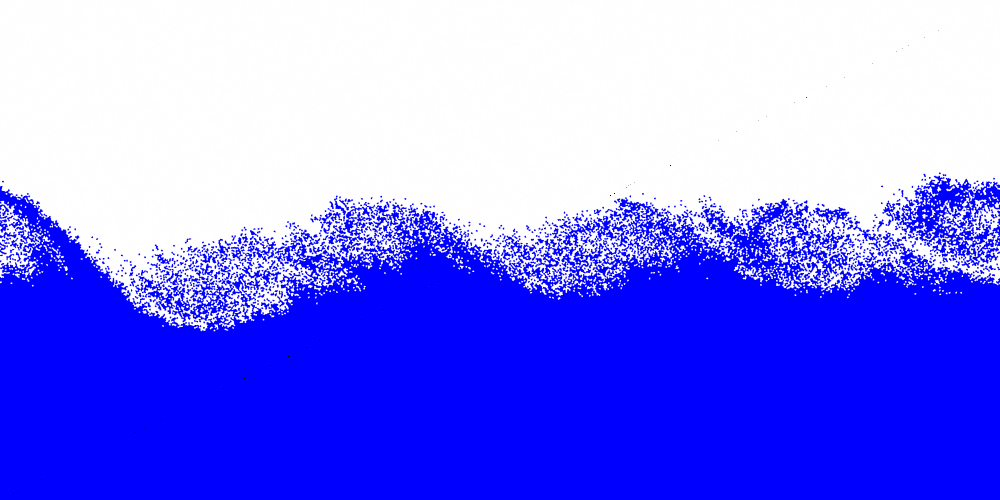}
  \caption{ Numerical illustration of the Kelvin-Helmotz instability
    on a rectangle with periodic conditions (in the horizontal
    coordinate) involving a discontinuous initial speed.  The
    parameters chosen for this experiment are given in
    \S\ref{subsec:RT}.
    \label{fig:kelvin-helmoltz} }
\end{figure}

For this second testcase, the domain is the rectangle $\Omega = [0,2]
\times [-.5,.5]$ periodized in the first coordinate by making the
identification identification $(4,x_2) \sim (0,x_2)$ for $x_2\in
[-.5,.5]$. The initial speed $v_0$ is discontinuous at $x_2 = 0$: the
upper part of the domain has zero speed, and the bottom part has unit
speed:
$$ v_0(x_1,x_2) = \left\{\begin{aligned}
0 \hbox{ if } x_2 \geq 0 \\
1 \hbox{ if } x_2 < 0
\end{aligned}\right.$$
This speed profile corresponds to a stationnary but unstable solution
to Euler's equation. If the subdomains $(\omega_i)_{1\leq i\leq N}$
are computed following \S\ref{subsec:optimomega}, the perfect symmetry
under horizontal translations is lost, and in
Figure~\ref{fig:kelvin-helmoltz} we observe the formation of vortices
whose radius increases with time. This experiment involves
$N=300\,000$ particles, with parameters $\tau = 0.005$ and $\eps=
0.0025$, and $2\,000$ timesteps.

\subsection{Rayleigh-Taylor instability}\label{subsec:RT}
For this last testcase, the particles are assigned a density $\rho_i$,
and are subject to the force of the gravity $\rho_i G$, where $G = (0,
-10)$.  This changes the numerical scheme to
\begin{equation} \label{eq:updrt}
  \left\{
  \begin{aligned}
    &V^{k+1}_i = V^k_i + \tau\left(\frac{1}{\eps^2}(\bary(\Lag_i(M^k,\psi^k) - M^k_i) + \rho_i G\right)\\
    &M^{k+1}_i = M^k_i + \tau V^{k+1}_i
  \end{aligned}
  \right.
\end{equation}
The computational domain is the rectangle $\Omega = [-1,1] \times
[-3,3]$, and the initial distribution of particles is given by $C_i =
\bary(\omega_i)$, where the partition $(\omega_i)_{1\leq i\leq N}$ is
constructed according to \S\ref{subsec:optimomega}. The fluid is
composed of two phases, the heavy phase being on top of the light phase:
$$ \rho_i = \left\{
\begin{aligned}
  &3 &\hbox{ if } C_{i2} > \eta \cos(\pi C_{i1}) \\
  &1 &\hbox{ if } C_{i2} \leq  \eta \cos(\pi C_{i1}) 
\end{aligned}
\right.,
$$ where $\eta = 0.2$ in the experiment and where we denoted $C_{i1}$
and $C_{i2}$ the first and second coordinates of the point
$C_{i}$. Finally, we have set $N = 50\,000$, $\eps = 0.002$ and $\tau
= 0.001$ and we have run $2000$ timesteps. The computation takes less
than six hours on a single core of a regular laptop. Note that it does
not seem straighforward to adapt the techniques used in the proofs of
convergence presented here to this setting, where the force depends on
the density of the particle. Our purpose with this testcase is only to
show that the numerical scheme behaves reasonably well in more complex
situations.

\begin{figure}
  \includegraphics[width=.15\textwidth]{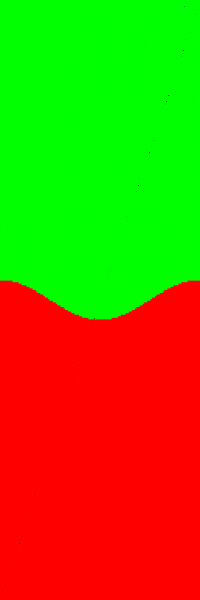} \hspace{.05cm} 
  \includegraphics[width=.15\textwidth]{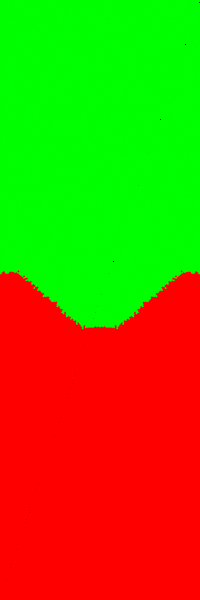} \hspace{.05cm}
  \includegraphics[width=.15\textwidth]{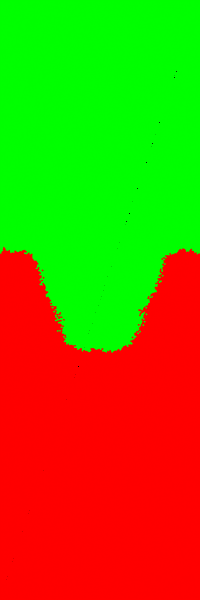} \hspace{.05cm}
  \includegraphics[width=.15\textwidth]{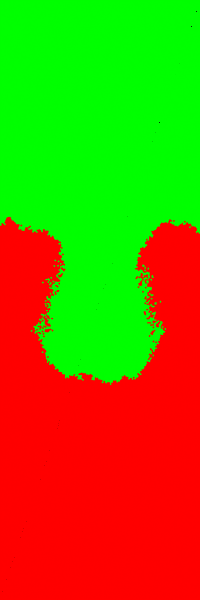} \hspace{.05cm}
  \includegraphics[width=.15\textwidth]{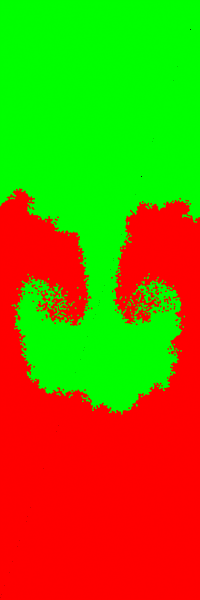} \hspace{.05cm}
  \includegraphics[width=.15\textwidth]{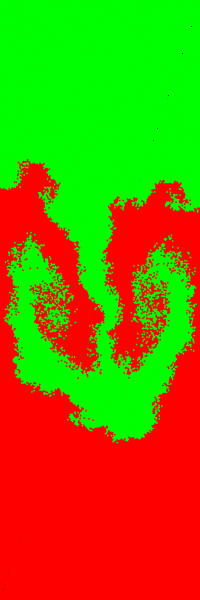} \\
  \includegraphics[width=.15\textwidth]{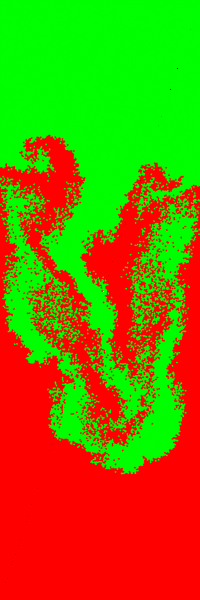} \hspace{.05cm} 
  \includegraphics[width=.15\textwidth]{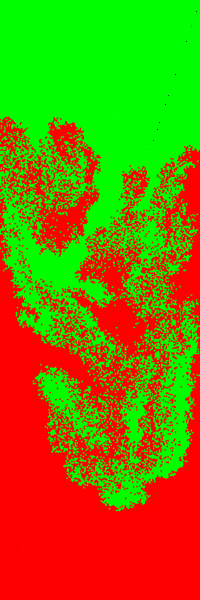} \hspace{.05cm} 
  \includegraphics[width=.15\textwidth]{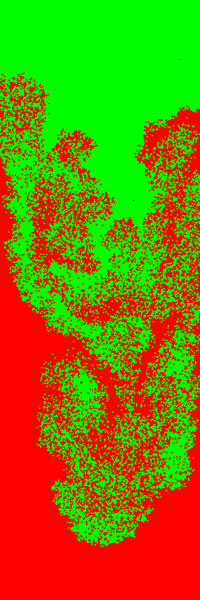} \hspace{.05cm} 
  \includegraphics[width=.15\textwidth]{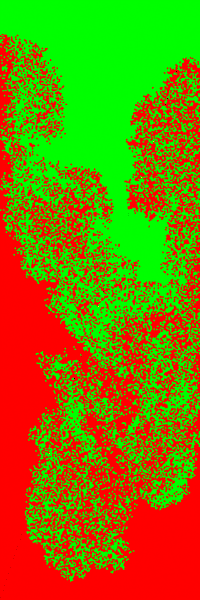} \hspace{.05cm} 
  \includegraphics[width=.15\textwidth]{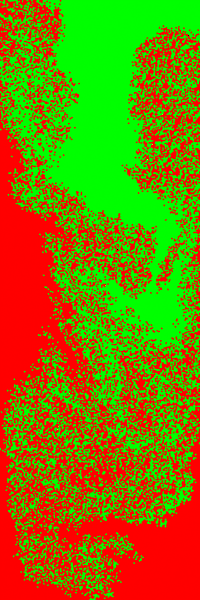}
  \caption{Numerical illustration of the Rayleigh-Taylor instability
    occuring when a heavy fluid (in green) is placed over a lighter
    fluid (in red). The parameters chosen for this experiment are given
    in \S\ref{subsec:RT}.
    \label{fig:rayleigh-taylor} }
\end{figure}

\bibliographystyle{abbrv}

\end{document}